\documentclass[a4paper,11pt]{article}

\usepackage[backend=biber,
    style=alphabetic,
    citestyle=alphabetic,
    doi=true,
    url=false,
    isbn=false,
    maxnames=999,
    abbreviate=true,
    giveninits=true]{biblatex}

\usepackage{csquotes}
\usepackage[tbtags]{amsmath}
\usepackage{amsthm}
\usepackage{aligned-overset}
\usepackage{amssymb, mathtools}
\newcommand{\negphantom}[1]{\settowidth{\dimen0}{#1}\hspace*{-\dimen0}}
\usepackage[nopatch=footnote]{microtype}
\DeclareMathSizes{10.95}{10}{8}{6}
\usepackage{marginnote}

\usepackage{graphicx}
\usepackage{caption}
\usepackage{subcaption}
\usepackage{float}
\usepackage[dvipsnames]{xcolor}
\usepackage{enumitem}
\usepackage[mathcal]{euscript}
\usepackage{tikz}
\usepackage{orcidlink}
\usepackage[a4paper%
  ,hcentering%
  ,textwidth=16cm%
  ,textheight=22.5cm%
  ,headheight=14.5pt%
  ,footskip=1.5cm%
  ,marginparwidth=.7in%
]{geometry}
\usepackage{hyperref}

\usepackage{zref-clever}
\zcRefTypeSetup{assumption}{
Name-sg = Assumption ,
name-sg = assumption ,
Name-pl = Assumptions ,
name-pl = assumptions ,
}
\AddToHook{env/thm/begin}{\zcsetup{countertype = {thm = theorem}}}
\AddToHook{env/propo/begin}{\zcsetup{countertype = {thm = proposition}}}
\AddToHook{env/lem/begin}{\zcsetup{countertype = {thm = lemma}}}
\AddToHook{env/cor/begin}{\zcsetup{countertype = {thm = corollary}}}
\AddToHook{env/ass/begin}{\zcsetup{countertype = {thm = assumption}}}
\AddToHook{env/definition/begin}{\zcsetup{countertype = {thm = definition}}}
\AddToHook{env/remark/begin}{\zcsetup{countertype = {thm = remark}}}
\zcsetup{cap = true}

\theoremstyle{plain}
\newtheorem{thm}{Theorem}[section]
\newtheorem{propo}[thm]{Proposition}
\newtheorem{lem}[thm]{Lemma}

\newtheorem{ass}[thm]{Assumption}
\newtheorem{definition}[thm]{Definition}
\newtheoremstyle{remark}
{3pt}
{3pt}
{}
{}
{\bfseries}
{}
{.5em}
{}
\theoremstyle{remark}
\newtheorem{remark}[thm]{Remark}

\newcommand{\sigmp}[1]{\;
  \begin{tikzpicture}[baseline=(content.base), inner sep=0pt]
    \node[inner sep=0pt] (content) {$\displaystyle #1$};
    \draw ([xshift=-0.1cm]content.west) -- ([xshift=-0.1cm, yshift=0.05cm]content.north west) -- ++(0.15,0);
    \draw  ([xshift=-0.1cm]content.west) -- ([xshift=-0.1cm,yshift=-0.05cm]content.south west) -- ++(0.15,0);
    \draw ([xshift=0.1cm]content.east) -- ([xshift=0.1cm, yshift=0.05cm]content.north east) -- ++(-0.15,0);
    \draw  ([xshift=0.1cm]content.east) -- ([xshift=0.1cm, yshift=-0.05cm]content.south east) -- ++(-0.15,0);
    \filldraw ([xshift=-0.1cm]content.north west) -- ([xshift=-0.05cm]content.west) --  ([xshift=-0.1cm]content.south west)  -- cycle;
    \filldraw ([xshift=0.1cm]content.north east) -- ([xshift=0.05cm]content.east) --  ([xshift=0.1cm]content.south east)  -- cycle;
    \node at ([xshift=0.2cm, yshift=0.03cm]content.north east) (C) {$\scriptscriptstyle \pm$};
    \vspace{-0.2cm}
  \end{tikzpicture}
  \hspace{-0.05cm}
}

\newcommand\bs[1]{\boldsymbol{#1}}
\newcommand{\R}{\mathbb{R}}
\newcommand{\diff}{\mathrm{\hspace{2pt}d}}
\newcommand{\X}{\bs{x}}
\newcommand{\dS}{\diff \mathcal{H}^{2}}
\newcommand{\bn}{\bs{n}}
\newcommand\norm[1]{\left\lVert #1 \right\rVert}
\newcommand\snorm[1]{\big\lVert #1 \big\rVert}

\newcommand\ssqrt[1]{\sqrt{\phantom{\mbox{$\cpm$}}\negphantom{\mbox{$\cpm$}}\smash{#1}}}

\newcommand{\bv}{\bs{v}}
\newcommand{\cpm}{c^{\pm}}
\newcommand{\psii}{\psi_i}
\newcommand{\psie}{\psi_e}
\newcommand{\lamMa}{\bs{\Lambda}(\bd)}
\newcommand{\varMa}{\bs{\mathcal{E}}(\bd)}
\newcommand{\varMadelta}{\bs{\mathcal{E}}(\bd_\delta)}
\newcommand{\bd}{\bs{d}}
\newcommand{\Mcont}{R}

\newcommand{\J}{\mathcal{J}}
\newcommand{\redJ}{\mathcal{J}_{\text{red}}}
\newcommand{\CXad}{\mathcal{X}_{\text{ad}}}
\newcommand{\CXR}{\mathcal{X}_R}
\newcommand{\CY}{\mathcal{Y}}
\newcommand{\CZ}{\mathcal{Z}}
\newcommand{\Sol}{\mathcal{S}}
\newcommand{\alpv}{\alpha_{\bv}}
\newcommand{\alpcpm}{\alpha_{\pm}}
\newcommand{\alpp}{\alpha_{+}}
\newcommand{\alpm}{\alpha_{-}}
\newcommand{\alpxi}{\alpha_{\xi}}
\newcommand{\lbxi}{\xi_{l}}
\newcommand{\ubxi}{\xi_{u}}

\newcommand{\lv}{\bs{\omega}}
\newcommand{\lp}{\pi}
\newcommand{\lcpm}{\gamma^{\pm}}
\newcommand{\lc}{\gamma}
\newcommand{\lpsi}{\Upsilon}
\newcommand{\lxi}{h}

\newcommand{\bbv}{\overline{\bs{v}}}
\newcommand{\bc}{\bar{c}}
\newcommand{\bcpm}{\bar{c}^{\pm}} 
\newcommand{\bpsi}{\bar{\psi}}
\newcommand{\bvh}{\bv_{\lxi}}
\newcommand{\ch}{c_{\lxi}}
\newcommand{\cpmh}{\cpm_{\lxi}}
\newcommand{\psih}{\psi_{\lxi}}
\newcommand{\CL}{\mathcal{L}}

\newcommand{\av}{\bs{\mathrm{u}}}
\newcommand{\ap}{\mathrm{q}}
\newcommand{\acpm}{\mathrm{r^{\pm}}}

\newcommand{\apsi}{\mathrm{s}}

\makeatletter
\newif\ifsubeqgroupstopped
\AtBeginDocument{\let\orig@theequation\theequation\let\orig@theHequation\theHequation}
\newcommand{\newsubeqgroup}[1]{\ifmeasuring@\else
  \setcounter{equation}{\value{parentequation}}%
  \global\let\theequation\orig@theequation\global\let\theHequation\orig@theHequation
  \refstepcounter{equation}\ltx@label{#1}%
  \protected@xdef\theparentequation{\theequation}%
  \protected@xdef\theHparentequation{\theHequation}%
  \setcounter{parentequation}{\value{equation}}%
  \setcounter{equation}{0}%
  \gdef\theequation{\theparentequation\alph{equation}}%
  \gdef\theHequation{\theHparentequation\alph{equation}}\fi}
\newcommand{\stopsubeqgroup}{\ifmeasuring@\else
  \setcounter{equation}{\value{parentequation}}%
  \global\let\theequation\orig@theequation\global\let\theHequation\orig@theHequation
  \global\subeqgroupstoppedtrue\fi}
\AtEndEnvironment{subequations}{\ifsubeqgroupstopped
  \setcounter{parentequation}{\value{equation}}\global\subeqgroupstoppedfalse\fi}
\makeatother

\begin{document}

\begin{center}
		
    {\Large \bf Optimal control for an anisotropic\\[0.2cm]Stokes--Nernst--Planck--Poisson system}

    \vskip0.7cm
    
    {\large\textsc{Giulia Cavalleri\orcidlink{0009-0006-4154-9659}$^1$}} \\
    {\normalsize e-mail: \texttt{cavalleri@wias-berlin.de}} \\
    \vskip0.35cm

    {\large\textsc{Luisa Plato\orcidlink{0000-0003-3456-1570}$^1$}} \\
    {\normalsize e-mail: \texttt{luisa.plato@wias-berlin.de}} \\
    \vskip0.35cm

    {\footnotesize $^1$Weierstrass Institute \\Anton-Wilhelm-Amo-Str. 39 \\ 10117 Berlin \\ Germany}
    \vskip0.7cm
  
\end{center}

\begin{abstract}
    \noindent We study an optimal control problem aimed at enhancing mixing in anisotropic electrolytes described by an instationary Stokes--Nernst--Planck--Poisson system in three space dimensions, where the anisotropy is encoded by a given director field in the diffusion and permittivity matrices. 
The control is the external electric potential entering through the Robin boundary condition of the Poisson equation. 
The cost functional combines tracking-type terms for the charge concentrations and an $L^2$-regularization of the control with a negative vorticity term that favors rotational flow. 
Assuming that the diffusion and permittivity matrices coincide, we first prove existence and uniqueness of weak solutions with nonnegative, essentially bounded concentrations, together with Lipschitz continuous dependence on the control. 
Since, through the negative vorticity term, the cost functional is no longer weakly lower semicontinuous, existence of an optimal control relies on strong convergence of the velocity gradients. 
Finally, we prove Fréchet differentiability of the control-to-state operator, introduce the adjoint system, and derive first-order necessary optimality conditions in the form of a variational inequality. 
\end{abstract}

\section{Introduction}

A \textit{lab-on-a-chip} is a device of a few square centimeters that handles fluid volumes of microliters or less and integrates laboratory processes such as sample transport, mixing of reagents, reaction, and detection on a single chip \cite{abgrall_2007}.
The concept dates back to 1979, when \cite{Terry_etal_79} used microelectronic techniques to build a microscale gas chromatograph.
Nowadays, applications of lab-on-a-chip range from DNA analysis and pathogen detection to cell manipulations \cite{abgrall_2007}, yet challenges such as specific sample preparation, device fabrication with larger solvent compatibility, and fluid handling without external pumps remain, making lab-on-chip an active field of research \cite{battat_2022}.
Although these devices can be built from different materials and rely on different technologies, they share common features, most notably their small size (down to the micro- and nano-scale). The advantages of portability and cost reduction are obvious, but miniaturization also causes some technical difficulties in the manipulation of fluids and particles. In this paper, we focus on the problem of mixing. At this scale, the Reynolds number is small, meaning that viscous forces dominate inertial ones. Thus, the flow is laminar rather than turbulent, and the mechanical stirring used in macroscopic devices is impractical here. Molecular diffusion, which is a slow process, prevails, and pressure-driven flow is not efficient, again due to the small diameter of the channels (see, e.g., \cite{zhao_yang_2012}). Achieving a good mixing of chemical reagents is therefore not an easy task.
A promising technique to overcome these challenges consists in dissolving ions in the fluid, thus obtaining an electrolyte, i.e., a fluid containing electrically charged particles that are free to move. The flow can then be manipulated by an external electric potential, which is generated by electrodes placed on the boundary and acts on the charged particles, achieving the desired mixing \cite{zhao_yang_2012}. This is particularly convenient since, in a microfluidic device, it is difficult to act on the fluid directly, for instance through a force distributed in its interior, whereas the potential applied at the electrodes can be easily tuned.
In this setting, three major effects govern the movement of the charges: diffusion, transport by the fluid, and electromigration (i.e., the movement induced by the electric field). In turn, space charges exert an electric body force on the fluid. 
Moreover, the electric field, the charges, and (consequently) the velocity field are influenced by the underlying anisotropy of the medium, that is, by the presence of a preferred direction, which may vary from point to point.

\paragraph{The state system.}
We consider a smooth and bounded domain $\Omega \subseteq \R^3$ with outer unit normal $\bn$ and boundary $\Gamma := \partial \Omega$. We assume that it is occupied by an incompressible fluid with velocity $\bv$, in which some positive ions with concentration $c^+$ and negative ions with concentration  $c^-$ are dissolved. An external electric potential $\xi$ is applied at the boundary. It may vary in space and time, influencing the bulk electric potential $\psi$ (and, in turn, the movement of the charges and of the fluid). Finally, the preferred direction of the medium at each point is given by the assigned director 
 $\bd: \overline{\Omega} \to \R^3$.
 An attempt to model the physical phenomena described above leads us to the instationary Stokes--Nernst--Planck--Poisson (SNPP) system (see, e.g., \cite[Chap.~3.4]{probstein}):
\begin{subequations}\label{eq:snpp}
	\begin{align}
			\kappa \, \partial_t \bv - \Delta \bv + \nabla p + (c^+ - c^-) \nabla \psi &= 0 &&\text{ in } Q \coloneqq (0,T) \times \Omega, \label{eq:s}\\
		\nabla \cdot \bv &= 0 &&\text{ in } Q, \label{eq:incomp}\\
		\partial_t c^{\pm} + \nabla \cdot (c^{\pm} \bv)
			- \nabla \cdot \left( \lamMa (
			\nabla c^{\pm} \pm c^{\pm} \nabla \psi) \right) &= 0 &&\text{ in } Q, \label{eq:np}\\
		- \nabla \cdot \left( \varMa
			\nabla \psi \right) - (c^+ - c^-)
			&= 0 &&\text{ in } Q, \label{eq:p}
\intertext{coupled with the boundary conditions\newsubeqgroup{eq:bc}}
		\bv &= 0 &&\text{ on } \Sigma \coloneqq (0,T) \times \Gamma, \label{eq:bn_s}\\
		\lamMa\left( \nabla \cpm \pm \cpm \nabla \psi \right) \cdot \bs{n} &= 0 &&\text{ on } \Sigma, \label{eq:bc_np}\\
		\varMa \nabla \psi \cdot \bn + \tau \psi &= \xi &&\text{ on } \Sigma, \label{eq:bc_p}
\intertext{and the initial conditions\stopsubeqgroup}
		\cpm(0) = \cpm_0, \qquad \bv(0) &= \bv_0 &&\text{ in } \Omega. \label{eq:ic}
    \end{align}
\end{subequations}
The velocity is described by the instationary Stokes system \eqref{eq:s}--\eqref{eq:incomp} with no-slip boundary condition \eqref{eq:bn_s}. Since, as pointed out above, the flow is dominated by viscosity (i.e., the Reynolds number is small), the convective term of the Navier--Stokes equations can be neglected, and the fluid is described by the Stokes equations.  Whether the time derivative has to be retained depends on the time scale of the forcing compared with the viscous diffusion time \cite[Sec.~6.1]{KimKarrila1991}. In our setting, the relevant time scale is set by the oscillation of the external electric potential $\xi$, which we assume to be rather short. Therefore, we keep the time derivative, i.e., we impose $\kappa > 0$. Besides this physical motivation, the time derivative is also important from a mathematical point of view, since it provides additional regularity of the velocity field needed in the analysis below, see Remark~\ref{rem:kappa}.

The charge densities evolve according to the Nernst--Planck equations \eqref{eq:np} including an anisotropic diffusion term $- \nabla \cdot(\lamMa \nabla \cpm)$ as well as two transport terms, one due to the velocity field $\bv$ and one due to the electric field $- \nabla \psi$. We impose the no-flux boundary condition \eqref{eq:bc_np} for the charge densities, meaning that we assume that the charges cannot cross the boundary of the domain. Here and throughout the paper, the symbol $\pm$ means an enumeration, so that \eqref{eq:np} denotes two equations: one for the positive charges $c^+$ (with a plus in front of the $\nabla \psi$ term) and one for the negative charges $c^-$ (with a minus in front of $\nabla \psi$). 

Finally, the electric potential is given by the Poisson equation \eqref{eq:p}.
As in \cite{bothe_fischer_saal_2014}, we choose Robin boundary conditions with $\tau>0$ \eqref{eq:bc_p} to model the electric double layer, which typically forms at an electrolyte-solid interface \cite[Chap.~7]{newman_thomas_alyea_2004}. More explicitly, when the electrolyte distributed in $\Omega$ is surrounded by a solid wall, its surface $\Gamma$ usually becomes charged: ions of opposite sign attach to the charged wall and form a thin layer that separates the wall from the rest of the electrolyte. This layer behaves like a capacitor, that is, the charge it stores is proportional to the difference between the external potential applied at the wall and the potential $\psi$ in the electrolyte. Equating this charge with the electric field flux $\varMa \nabla \psi \cdot \boldsymbol{n}$ yields \eqref{eq:bc_p}, where $\tau$ is the local capacity and $\xi / \tau$ is the potential applied at the wall.
The matrices $\lamMa = I + \lambda \, \bd \otimes \bd$ and $\varMa = I + \varepsilon \, \bd \otimes \bd$ with $\lambda, \varepsilon > 0$ appearing in the Nernst--Planck and Poisson equations are the anisotropy matrices. They have been chosen in this form, referring to \cite{calderer_2016}, with the idea that the present model could be extended to describe nematic liquid crystals, where $\bd$ is also time-dependent and ruled by an additional PDE (see \cite{banas_lasarzik_prohl_2021}). However, in the following, we will impose $\lamMa = \varMa$ for mathematical reasons (see \zcref{remark:lam_var}) as already done in the literature (see \cite{banas_lasarzik_prohl_2021} and \cite{FeireislRocca}). 

\paragraph{Optimal mixing.}
Having the design of lab-on-a-chip devices in mind, we would like to maximize the mixing of the electrolyte described by the initial-boundary value problem \eqref{eq:snpp}--\eqref{eq:ic} by controlling the external electric potential $\xi$ applied at the boundary. 
To achieve this goal, we employ two complementary strategies.
On the one hand, we aim to maximize the vorticity of the velocity field $\bv$,  which increases the mixing by enhancing the local rotation of the fluid, as has already been done in \cite{Hu_Wu_2018}.
On the other hand, we want to minimize the variance of the charge concentrations $\cpm$, which is a measure of the spatial homogeneity of the concentrations.
Mathematically, we introduce the \textit{cost functional}
\begin{equation}
    \begin{split}\label{eq:obj}
         \J(\xi, (\bv, \cpm, \psi)) \coloneqq &\int_0^T \!\!\int_\Omega \Big(- \frac{\alpv}{2}|\nabla \times \bv|^2 + \frac{\alpp}{2}|c^+ - c^+_Q|^2 + \frac{\alpm}{2}|c^- - c^-_Q|^2  \Big)\diff \X + \frac{\alpxi}{2} \int_\Gamma |\xi|^2 \dS \! \diff t, 
    \end{split}
\end{equation}
where $\alpv, \alpcpm, \alpxi$ are nonnegative weights that cannot vanish simultaneously and $c^{\pm}_Q$ are given target functions. 
The first term in the cost functional \eqref{eq:obj} is minus the $L^2(Q)$-norm of the vorticity: due to the negative sign in front of it, minimizing it means maximizing the vorticity.
The second and third terms in the cost functional are of standard tracking type, where the $c^{\pm}_Q$ are given evolution targets. 
For our analysis, a relevant case is the constant function $c^{\pm}_Q$ equal to the spatial mean of the initial data: an electrolyte is well-mixed precisely when each species is spread uniformly, i.e., spatially constant and equal to its mean, and, as we will see (cf. \eqref{eq:const_mean}), the means are conserved throughout the evolution, so $\frac{1}{|\Omega|} \int_\Omega\cpm_0 \diff \X$ is the correct target at every time. 
As already anticipated, with that choice of target function $c^{\pm}_Q$, the tracking type terms for the charges correspond to the variance of the concentrations.
The last term in the cost functional \eqref{eq:obj} is an $L^2(\Sigma)$-regularization for the control.
We aim to study the following optimal control problem
\begin{gather}
    \textit{Minimize the cost functional } \J \textit{ subject to the state system \eqref{eq:snpp}--\eqref{eq:ic} and to the control constraint} \nonumber\\
    \xi \in \CXad \coloneqq\{\xi \in L^{\infty}(\Sigma) \, \big| \; \lbxi \leq \xi \leq \ubxi \text{ a.e. in } \Sigma \}, \label{eq:adm_set}
\end{gather}
where $\lbxi, \ubxi$ are given $L^{\infty}$-functions. After proving the existence of minimizers, our goal is to characterize them through first-order necessary conditions, which open the route to numerical simulations. 

\begin{remark}\label{rem:kappa}
    Since cost functionals of this type are not weakly lower semicontinuous with respect to the curl of the velocity, weak convergence of a minimizing sequence does not suffice to pass to the limit in the cost functional in the existence proof for an optimal control. 
    Instead, we require strong convergence of the first-order derivatives of the velocity in $L^2(0,T;L^2(\Omega))$ along a subsequence. 
    We obtain this from an Aubin--Lions compactness argument using the maximal $L^2$-regularity of the instationary Stokes equations, which gives a uniform bound for the velocity in $W^{1,2}(0,T;L^2(\Omega)) \cap L^2(0,T;W^{2,2}(\Omega))$, which compactly embeds into $L^2(0,T;W^{1,2}(\Omega))$. It is precisely here that the assumption $\kappa > 0$ is essential.
\end{remark}

\begin{remark}
    The admissible control set $\CXad$ is a convex and bounded subset of $L^{\infty}(\Sigma)$. 
    In particular, there exists a certain open ball
    \begin{equation*}
         \CXR \coloneqq\{\xi \in L^{\infty}(\Sigma)\, \big| \; \|\xi\|_{L^{\infty}(\Sigma)} < \Mcont \}
    \end{equation*}
    such that $\CXad \subseteq \CXR$. 
    Moreover, $\CXad$ is weakly-$\ast$ closed. 
    The proof is simple. Take a sequence $\{\xi_n\}_n \subseteq \CXad$ such that $\xi_n \overset{\ast}{\rightharpoonup} \xi$ in $L^{\infty}(\Sigma) = (L^1(\Sigma))^*$. 
    Then, $\xi_n\le\ubxi$, and we have
    \begin{equation*}
        \int_0^{T} \int_{\Gamma} (\ubxi-\xi_n)\varphi \dS \diff t \to \int_0^{T} \int_{\Gamma} (\ubxi-\xi)\varphi \dS \diff t 
    \end{equation*}
    for every nonnegative $\varphi \in L^1(\Sigma)$. Since the sequence of integrals is nonnegative, the limit is, in turn, nonnegative. 
    Thus, $\ubxi - \xi \geq 0$. The lower bound is analogous, and we have proved that $\xi \in \CXad$.
\end{remark}

\paragraph{Related literature.} 
The introduction of the coupled Navier--Stokes--Nernst--Planck--Poisson system goes back to the work of Rubinstein in 1990 \cite{rubinstein}, and the mathematical literature nowadays counts several contributions on the related analytical results (see \cite{jerome_2002,schmuck_2009,bothe_fischer_saal_2014,constantin_ignatova_2019,constantin_ignatova_lee_2021,lee_2022,homberg_lasarzik_plato_2025} and the references cited therein but also \cite{lasarzik_plato_2025} for the stationary Stokes--Nernst--Planck--Poisson system). Regarding the optimal control problem, tracking-type cost functionals are common, but vorticity-based cost functionals have been less studied. In \cite{abergel_temam_1990}, the $L^2$-norm of the vorticity is minimized in order to reduce turbulence; see also \cite{kasumba_kunisch_2012}.
The maximization of vorticity was considered in the context of shape optimization in \cite{simon_notsu_2022_stokes, simon_notsu_2022_ns} and in the context of mixing in \cite{Hu_Wu_2018}.
Optimal control problems for nonlinear complex fluids have been studied intensively over the last decade, and we do not attempt to 
give a comprehensive review here. However, we mention some relevant results for the context of the present paper.
For a general and easy introduction to the analysis of optimal control problems, we refer to \cite{Troltzsch_2010}.
Navier--Stokes--Cahn--Hilliard systems of fourth and sixth order have been studied, respectively, in \cite{hintermuller_wegner_2014, frigeri_rocca_sprekels_2016} and in the recent \cite{Signori_Wu_26}, where the authors prove the existence of optimal controls and first-order optimality conditions. Notice that the extension from a single fluid to a two-phase flow of immiscible fluids, which can be described by coupling the system with a Cahn--Hilliard equation for a diffuse interface, is a natural step from the point of view of applications.
The optimal control of chemotaxis--Navier--Stokes models has been studied in \cite{lopez-rios_villamizar-roa_2021}, 
where the Navier--Stokes equations are coupled to certain transport diffusion equations for the chemotactic species, 
which is in some sense similar to transport diffusion equations for the charges in the present paper, but there distributed controls instead of boundary controls are used.
An optimal control problem in the context of semiconductors was studied in \cite{hinze_pinnau_2002}, 
where the authors study the optimal control of the doping profile in a drift-diffusion model for semiconductors 
using a stationary Nernst--Planck equation for the charge density and a Poisson equation for the potential.
The numerical optimal mixing of a micromixer was studied in \cite{ji_deng_2016} based on a stationary Navier--Stokes--Poisson system,
where the external potential is also used as a control.
However, to the best of our knowledge, a rigorous analysis of optimal control problems for Stokes--Nernst--Planck--Poisson systems
has not yet been carried out, neither for the isotropic nor for the anisotropic setting.
The present paper aims to close this gap.

\paragraph{Main results and outline of the paper.} 
The paper follows the classical route in optimal control theory, showing the well-posedness of the state system and the differentiability of the control-to-state operator to then derive first-order optimality conditions.
The main difficulties lie in the well-posedness of the state system and the differentiability of the control-to-state operator.
Our main contributions are the following.
\begin{itemize}
    \item \emph{Well-posedness for rough boundary controls.}
    We prove the existence and uniqueness of global-in-time weak solutions to the anisotropic SNPP system \eqref{eq:snpp}--\eqref{eq:ic} 
    for controls $\xi \in L^\infty(\Sigma)$ under the structural assumption $\lambda = \varepsilon$ (\zcref{thm:weak_existence} and \zcref{thm:uniqueness}).
    The concentrations are nonnegative and essentially bounded, uniformly with respect to the admissible controls,
    and the solution depends Lipschitz continuously on the control and on the initial data.
    The existence results for similar (Navier--)Stokes--Nernst--Planck systems in the available literature exploit a more regular external potential \cite{homberg_lasarzik_plato_2025, bothe_fischer_saal_2014, Fischer_Saal_17},
    whereas the low regularity of $\xi$ is natural in the presence of box constraints.
    This requires new a~priori estimates, which rely on $W^{1,p}$-regularity for elliptic problems with Robin boundary conditions \cite{amrouche_conca_ghosh_ghosh_20},
    yielding $\psi \in L^\infty(0,T;W^{1,6}(\Omega))$,
    and on a Moser--Alikakos iteration \cite{alikakos_79} for the $L^\infty$-bounds of the concentrations.
    \item \emph{Existence of optimal controls for a cost functional that is not weakly lower semicontinuous.}
    Due to the negative vorticity term, weak convergence of solutions is not sufficient to pass to the limit in the cost functional. 
    Therefore, we need strong convergence in the velocity gradients.
    We achieve this convergence by exploiting the maximal $L^2$-regularity of the instationary Stokes equation,
    which yields strong convergence of the velocity gradients along minimizing sequences.
    \item \emph{Differentiability and first-order optimality conditions.}
    We prove the Fr\'echet differentiability of the control-to-state operator on an open ball in $L^\infty(\Sigma)$,
    establish the well-posedness of the adjoint system,
    and derive first-order necessary optimality conditions in the form of a variational inequality (\zcref{thm:opt_cond}).
\end{itemize}

The paper is organized as follows. In Section~\ref{sec:notation_preliminaries} we introduce some notation and collect some useful inequalities used later on in the paper.
In Section~\ref{sec:well_posedness} we state and prove the well-posedness of the state system, i.e., the existence and uniqueness of a global-in-time weak solution and the continuous dependence on the initial and boundary data.
In Section~\ref{sec:optimal_control} we prove the existence of an optimal control.
In Section~\ref{sec:differentiability} we prove the differentiability of the control-to-state operator, and in Section~\ref{sec:first_order} we conclude the work by deriving first-order optimality conditions.

\section{Notation and preliminaries}\label{sec:notation_preliminaries}
\paragraph{Constants.}
Throughout, $C$ denotes a generic large positive constant and $\delta$ a generic small positive constant, whose value may change from line to line and which depends only on $\Omega$, $T$, $\lambda$, $\varepsilon$, $\tau$, $\kappa$, $\bd$. Dependence on further quantities is always explicitly stated and, whenever we need to keep track of a parameter, indicated by subscripts (as in $C_{p}$, $C_{\delta}$, $C_{R}$ or $C_{R,\delta}$).

\paragraph{Scalar and tensor products.}
Throughout the paper, vectors and matrices are written in bold, and 
\begin{equation*}
    \bs{a}\cdot\bs{b} = \sum_{i=1}^3 a_i b_i, \qquad \bs{a}\otimes\bs{b} =  \big( a_i b_j \big)_{i,j=1}^3, \qquad \bs{A}:\bs{B} = \sum_{i,j=1}^3 A_{ij} B_{ij}
\end{equation*}
denote, respectively, the scalar product of two vectors, their tensor product, and the Frobenius product of two matrices.

\paragraph{Function spaces.}
We employ the standard notation for the Lebesgue and Sobolev spaces. We use $W^{1,2}_{0}(\Omega)$ to denote the space of $W^{1,2}(\Omega)$-functions having vanishing trace at the boundary, and we set 
    \begin{equation*}
        \begin{split}
            W^{1,2}_{0,\sigma}(\Omega;\R^3) &\coloneqq \big\{\bs{\varphi}\in W^{1,2}_{0}(\Omega;\R^3)\ :\ \nabla\cdot\bs{\varphi}=0\big\},\\
        L^{2}_{\sigma}(\Omega;\R^3)& \coloneqq\overline{W^{1,2}_{0,\sigma}(\Omega;\R^3)}^{L^2(\Omega)}
        \end{split}
    \end{equation*}
for the vector-valued divergence-free spaces. 
Throughout the paper $\langle\cdot,\cdot\rangle$ denotes the duality between  $(W^{1,2}(\Omega))^*$ and $W^{1,2}(\Omega)$.
For $p\in(1,\infty)$ we write
    \begin{equation*}
        W^{-1/p,p}(\Gamma)\coloneqq (W^{1/p,p'}(\Gamma))^*,
    \end{equation*}
where $p'= p/(p-1)$.
In the following, to keep the notation as light as possible, we do not explicitly indicate the dimension of the target space in the function spaces, e.g., we write $\bv \in W^{1,2}(\Omega)$ instead of $\bv \in W^{1,2}(\Omega;\R^3)$. It will be clear from the context and from the use of the bold font whether scalar- or vector-valued functions are meant.

\paragraph{Useful inequalities.}

We collect here the main inequalities used throughout the paper. 

\begin{itemize}
    \item (Young's inequality) Let $a,b \geq 0$ be nonnegative real numbers and $p,p' \in (1, \infty) $ conjugate H\"older exponents. Then, for every $\eta > 0$,
    \begin{equation}\label{eq:young}
        a b \leq \eta a^p + C_{\eta} b^{p'}
    \end{equation}
    where $ \displaystyle C_{\eta} = \frac{1}{p' (p\eta)^{p'/p}}$.
    \item (Generalized Poincaré's inequality) Let $\Omega$ be a Lipschitz bounded domain. Then, there exists a positive constant $C_P>0$ depending on $\Omega$ such that
    \begin{equation} \label{eq:gen_Poincare_ineq}
        \| f \|_{L^2(\Omega)} \leq C_P ( \| \nabla f \|_{L^2(\Omega)} +  \| f \|_{L^2(\Gamma)}) \qquad \forall f \in W^{1,2}(\Omega).
    \end{equation}
    \item (Trace inequality) Let $\Omega$ be a bounded Lipschitz domain. Then, there exists a positive constant $C>0$ such that  
    \begin{align}
        \norm{f}_{L^2(\Gamma)} &\leq C \norm{f}_{L^2(\Omega)}^{1/2} \norm{f}_{W^{1,2}(\Omega)}^{1/2}, \label{eq:trace_L2}
    \end{align}
    for all $f \in W^{1,2}(\Omega)$.  
   It follows from \cite[Cor.~2, Section 1.4.7, p.~82]{mazya_11} applied---with the notation of the book---for $\mu = \mathcal{H}^2|_{\partial \Omega}$, $k=0, \ell = 1, p = 2$, $q = 2, \tau = 1/2$. Since this result is stated only for $f \in C^1(\Omega)$, it can then be extended to $f \in W^{1,2}(\Omega)$ by density and continuity of the trace operator. 

    \item (Special cases of Gagliardo--Nirenberg's inequality) Let $\Omega \subseteq \R^3$ be a Lipschitz bounded domain. Then, there exists a positive constant $C>0$ such that
    \begin{align}
        \norm{f}_{L^3(\Omega)} &\leq C \Big(\norm{f}_{L^2(\Omega)}^{1/2} \norm{\nabla f}_{L^2(\Omega)}^{1/2} + \norm{f}_{L^2(\Omega)} \Big) \qquad &\forall f \in W^{1,2}(\Omega), \label{eq:GN_L3} \\
        \norm{f}_{L^3(\Omega)} &\leq C \norm{f}_{L^2(\Omega)}^{1/2} \norm{\nabla f}_{L^2(\Omega)}^{1/2} \qquad &\forall f \in  W^{1,2}_0(\Omega), \label{eq:GN_L3_0}\\
        \norm{f}_{L^4(\Omega)} &\leq C \Big(\norm{f}_{L^2(\Omega)}^{1/4} \norm{\nabla f}_{L^2(\Omega)}^{3/4} + \norm{f}_{L^2(\Omega)} \Big) \qquad &\forall f \in W^{1,2}(\Omega), \label{eq:GN_L4} 
    \end{align}
\end{itemize}

We also prove the following result, which is an easy consequence of the Gagliardo--Nirenberg inequality. 

\begin{lem}\label{lem:gagliardo_nirenberg_L1}
    Let $\Omega \subseteq \R^3$ be a Lipschitz bounded domain. Then, there exists a positive constant $C>0$ such that
        \begin{align}\label{eq:gagliardo_nirenberg_L1}
        \norm{f}_{L^3(\Omega)}
        \leq C \Big(\norm{\nabla f}_{L^2(\Omega)}^{4/5} \norm{f}_{L^1(\Omega)}^{1/5}
        + \norm{f}_{L^1(\Omega)} \Big) \qquad \forall f \in W^{1,2}(\Omega).
    \end{align}
\end{lem}
\begin{proof}
    This is the Gagliardo--Nirenberg inequality
    $\norm{f}_{L^3(\Omega)} \leq C \norm{f}_{W^{1,2}(\Omega)}^{4/5} \norm{f}_{L^1(\Omega)}^{1/5}$,
    whose exponent is determined by $\frac{1}{3} = \frac{4}{5} \bigl( \frac{1}{2} - \frac{1}{3} \bigr) + \frac{1}{5}$,
    with the $L^2$-part of the $W^{1,2}$-norm eliminated.
    Indeed, the subadditivity of $t \mapsto t^{4/5}$, the interpolation
    $\norm{f}_{L^2(\Omega)} \leq \norm{f}_{L^1(\Omega)}^{1/4} \norm{f}_{L^3(\Omega)}^{3/4}$
    and Young's inequality with the exponents $(5/3, 5/2)$ give
    \begin{align*}
        C \norm{f}_{W^{1,2}(\Omega)}^{4/5} \norm{f}_{L^1(\Omega)}^{1/5}
        &\leq C \norm{\nabla f}_{L^2(\Omega)}^{4/5} \norm{f}_{L^1(\Omega)}^{1/5}
            + C \norm{f}_{L^1(\Omega)}^{2/5} \norm{f}_{L^3(\Omega)}^{3/5}\\
        &\leq C \norm{\nabla f}_{L^2(\Omega)}^{4/5} \norm{f}_{L^1(\Omega)}^{1/5}
            + \frac{1}{2} \norm{f}_{L^3(\Omega)}
            + C \norm{f}_{L^1(\Omega)} ,
    \end{align*}
    where $\frac{1}{2} \norm{f}_{L^3(\Omega)}$ is absorbed into the left-hand side of \eqref{eq:gagliardo_nirenberg_L1}.
\end{proof}

\paragraph{The anisotropy matrices.}
For a vector field $\bd \in W^{1,\infty}(\Omega)$ the matrices
\begin{equation*}
    \lamMa = I+\lambda\,\bd\otimes\bd, \qquad \varMa = I+\varepsilon\,\bd\otimes\bd
\end{equation*}
are symmetric and belong to $W^{1,\infty}(\Omega;\mathbb{R}^{3\times3})$. For
$\bs{a}\in\mathbb{R}^{3}$ we use the shorthand
\begin{equation*}
    |\bs{a}|^{2}_{\lamMa} \coloneqq\lamMa\,\bs{a}\cdot\bs{a} = |\bs{a}|^{2}+\lambda\,|\bd\cdot\bs{a}|^{2},
\end{equation*}
and analogously for $\varMa$. Since $\lambda,\varepsilon>0$ and
$\bd\in W^{1,\infty}(\Omega)$, we have
\begin{equation*}
    |\bs{a}|^{2}\le|\bs{a}|^{2}_{\lamMa} \le\big(1+\lambda\|\bd\|^{2}_{L^{\infty}(\Omega)}\big)|\bs{a}|^{2}
    \qquad \forall\bs{x}\in\Omega,\ \bs{a}\in\mathbb{R}^{3},
\end{equation*}
and the same bounds hold for $\varMa$ with $\varepsilon$ in place of $\lambda$.

\paragraph{Compact notation for the charges.}
All equations for the two ionic species are written simultaneously by means of the superscript
$\pm$: an equation containing $c^{\pm}$ and the signs $\pm$, $\mp$ stands for the two
equations obtained by reading either all upper or all lower signs. Moreover, sums over the two species
are abbreviated by
\begin{equation*}
\sigmp{f^{\pm}}\coloneqq f^{+}+f^{-}.
\end{equation*}

\paragraph{A regularity result for the electric potential.}     
Here we recall a well-known regularity result for the Poisson equation with Robin boundary conditions, which we write under the assumptions we need in our case, but that holds in a stronger sense. For further details, the interested reader can refer to \cite{amrouche_conca_ghosh_ghosh_20}.
\begin{propo}[$W^{1,p}$-integrability of the electric potential]\label{propo:gradient_integrability} 
    Let $\Omega$ be a bounded $C^1$-domain, $p\in(1,\infty)$ and $q\in[1,\infty]$. Assume $F \in L^q(0,T;L^p(\Omega))$, $g \in L^\infty(\Sigma)$, and $\bd \in W^{1,\infty}(\Omega)$. Then, the solution $\psi$ of the Poisson equation with the Robin boundary condition
    \begin{equation*}
        \left\{
        \begin{aligned}
            -\nabla\cdot\big(\varMa\nabla\psi\big) &= F && \text{ in }\Omega,\\
             \varMa\nabla\psi\cdot\boldsymbol{n}+\tau\psi&=g\ &&\text{ on }\Gamma,
        \end{aligned}
        \right.
    \end{equation*}
     has the regularity $L^q(0,T;W^{1,p}(\Omega))$ and satisfies
    \begin{equation*}
        \|\psi\|_{L^q(0,T;W^{1,p}(\Omega))}\le C_p\Big(\big\|F\big\|_{L^q(0,T;L^p(\Omega))} +\|g\|_{L^\infty(\Sigma)}\Big)
    \end{equation*}
    for a positive constant $C_p>0$ depending on $p$, $\Omega$, $\tau$ and $\boldsymbol{d}$ but not on $F$ and $g$.
\end{propo}

\begin{proof}
    For almost every $t\in(0,T)$, the function $\psi(t)$ solves the Robin problem studied in \cite{amrouche_conca_ghosh_ghosh_20} with $A=\varMa$, $\boldsymbol{f}=0$, $\alpha=\tau$. The assumptions of
    \cite[Thm.~2.1]{amrouche_conca_ghosh_ghosh_20} are satisfied:
    \begin{itemize}
        \item $\boldsymbol{d}\in W^{1,\infty}(\Omega) \subseteq C^0(\bar{\Omega})$ by Sobolev embedding; therefore $\varMa=I+\varepsilon\,\boldsymbol{d}\otimes\boldsymbol{d}$ is continuous on $\overline{\Omega}$ and, extending it to a continuous matrix field on $\R^3$, it belongs to $\mathrm{VMO}(\R^3)$;
        \item for a.e. $t$, $F(t,\cdot) \in L^p(\Omega)\hookrightarrow\big(W^{1,p'}(\Omega)\big)^*$;
        \item for a.e. $t$, $g(t,\cdot) \in L^\infty(\Gamma)\hookrightarrow W^{-1/p,p}(\Gamma)$, since
        \begin{equation*}
            \begin{split}
                \bigg|\int_\Gamma g(t,\cdot)\varphi\,\diff\mathcal{H}^2\bigg| 
                &\le \|g(t,\cdot)\|_{L^\infty(\Gamma)}\|\varphi\|_{L^1(\Gamma)}
                \le C\|g\|_{L^\infty(\Sigma)}\|\varphi\|_{W^{1/p,p'}(\Gamma)}
            \end{split}
        \end{equation*}
        for every $\varphi\in W^{1/p,p'}(\Gamma)$, by the boundedness of $\Gamma$.
    \end{itemize}
    Thus $\psi(t)\in W^{1,p}(\Omega)$ and the corresponding stationary estimate holds, with a constant depending only on $p$, $\Omega$, $\tau$ and $\boldsymbol{d}$. Raising the stationary estimate to the power $q$ and integrating in time (or taking the essential supremum if $q=\infty$) gives the claim.
\end{proof}

\paragraph{Splitting of the electric potential.}
As is classically done in the literature (see, e.g., \cite{Fischer_Saal_17}), we split the electric potential into an internal part $\psii$, generated by the charges, and an external part $\psie$, generated by the control $\xi$:
\begin{equation}\label{eq:splitting}
    (P_i)\left\{
    \begin{aligned}
        - \nabla \cdot (\varMa \nabla \psii) &= c^+ - c^- && \text{in } \Omega, \\
        \varMa \nabla \psii \cdot \bn + \tau \psii &= 0 && \text{on } \Gamma,
    \end{aligned}
    \right.
    \quad
    \qquad
    (P_e)\left\{
    \begin{aligned}
        - \nabla \cdot (\varMa \nabla \psie) &= 0 && \text{in } \Omega, \\
        \varMa \nabla \psie \cdot \bn + \tau \psie &= \xi && \text{on } \Gamma,
    \end{aligned}
    \right.
    \quad
\end{equation}
so that $\psi = \psii + \psie$. Notice that $(P_e)$ does not involve the unknowns $(\bv, \cpm)$: the external potential $\psie$ is determined by the control $\xi$ alone and can therefore be regarded as given data. The whole coupling between the electric potential and the charges is contained in $(P_i)$, which has homogeneous boundary conditions. Moreover, they satisfy the following estimates.
\begin{itemize}
    \item  For $\psii$ we have full elliptic regularity: if $c^+ - c^- \in L^p(\Omega)$, then, by \cite[Thm.~2.4.2.7]{grisvard_85}, $\psii \in W^{2,p}(\Omega)$ with
    \begin{equation}\label{eq:est_psii}
        \norm{\psii}_{W^{2,p}(\Omega)} \leq C \norm{c^+ - c^-}_{L^p(\Omega)}
    \end{equation}
    almost everywhere in time. We are going to use this estimate with $p=2$.
    \item   For $\psie$,  we apply \zcref{propo:gradient_integrability} with $F = 0$, $g = \xi$, and $q=+\infty$. Then, the following inequality holds
    \begin{equation}\label{eq:est_psie}
        \norm{\psie}_{L^{\infty}(0,T;W^{1,p}(\Omega))} \leq C \norm{\xi}_{L^{\infty}(\Sigma)} \leq C R = \tilde{C}
    \end{equation}
    for a constant $\tilde{C}>0$ that depends on $R$ (and $p$). We are going to apply this estimate with $p=6$.
\end{itemize}
Then, $\psi = \psii + \psie$ has the lower of the two regularities. 

\section{Well-posedness of the state system}\label{sec:well_posedness}
Before stating the well-posedness results, we need to impose some assumptions on the data, which we collect in the following and that will hold throughout the paper.
\begin{ass}\label{ass:strong}
    We impose the following conditions.
    \begin{enumerate}[label = \textup{(A\arabic*)}]
        \item \textup{(Domain.)} 
        $\Omega \subset \R^3$ is a bounded domain with $C^2$-boundary $\Gamma \coloneqq\partial \Omega$.  \label{ass:domain}
        \item \textup{(Anisotropy.)}
        The director field $\bd \in W^{1,\infty}(\Omega)$.
        The matrices $\lamMa, \varMa \in W^{1,\infty}(\Omega;\R^{3 \times 3})$ are given by $\lamMa = I + \lambda \bd \otimes \bd$ and $\varMa = I + \varepsilon \bd \otimes \bd$ for some $\lambda, \varepsilon > 0$.  \label{ass:anisotropy}
        \item \textup{(Structural condition.)} 
        The diffusion and permittivity matrices coincide, $\lamMa = \varMa$, corresponding to $\varepsilon = \lambda$. \label{ass:struct_cond}
        \item \textup{(Robin parameter.)} 
         $\tau > 0$. \label{ass:robin_cond}
        \item \textup{(Initial data.)} 
        $c^{\pm}_0 \in L^\infty(\Omega)$ with $0 \leq c^\pm_0 \leq M^{\pm}$ almost everywhere in $\Omega$, $\bv_0 \in W^{1,2}_{0, \sigma}(\Omega)$. \label{ass:initial_data}
        \item \textup{(Cost functional and admissible set.)} 
        The weights $\alpv, \alpcpm, \alpxi$ are nonnegative and cannot vanish simultaneously, and $c^{\pm}_Q \in L^2(Q)$ are given target functions. The functions $\lbxi$ and $\ubxi$ in \eqref{eq:adm_set} belong to $L^{\infty}(\Sigma)$ and satisfy the compatibility condition $\lbxi < \ubxi$ a.e. in $\Sigma$.
        \label{ass:cost_funct}
    \end{enumerate}
\end{ass}

\begin{remark}\label{remark:lam_var}
    The difference between the anisotropy matrices $\varMa$ and $\lamMa$ is crucial for capturing the distinct physical properties of electric conductivity and electric permittivity (see \cite{calderer_2016}). However, the technical assumption $\varepsilon = \lambda$ in \ref{ass:struct_cond} is needed in the existence proof to rewrite the term $ \mp \nabla \cdot (\lamMa \nabla \psi)$ appearing in \eqref{eq:np} by means of the Poisson equation \eqref{eq:p} (see Steps \ref{enum:energy_est}--\ref{enum:grad_est_charges} in \zcref{sec:weak_exist}), which is crucial in our estimates.  From now on, we will use only $\varMa$, dropping the notation $\lamMa$.
\end{remark}

\begin{remark}
    In \ref{ass:initial_data} we require that $\bv_0 \in W^{1,2}_{0, \sigma}(\Omega)$ for the initial velocity to get the full maximal regularity of the nonstationary Stokes equation \cite[Section IV, Thm.~2.7.1]{sohr_01}.
\end{remark}

\noindent We start by defining our notion of weak solutions.

\begin{definition}[Weak solution]\label{definition:weak_sol}
    A tuple $(\bv, \cpm, \psi)$ is called a \emph{weak solution} of the SNPP initial-boundary value problem \eqref{eq:snpp}--\eqref{eq:ic} if
    \begin{align*}
        \bv &\in W^{1,2}(0,T;L^2(\Omega)) \cap L^2(0,T;W^{1,2}_{0,\sigma}(\Omega) \cap W^{2,2}(\Omega)),\\
        \cpm &\in W^{1,2}(0,T;(W^{1,2}(\Omega))^*) \cap L^\infty(0,T;L^2(\Omega)) \cap L^2(0,T;W^{1,2}(\Omega)),\\
        \psi &\in L^\infty(0,T;W^{1,6}(\Omega)),
    \end{align*}
    where the initial conditions \eqref{eq:ic} hold in $L^2(\Omega)$, and the system \eqref{eq:snpp}--\eqref{eq:bc}
    is satisfied in the following variational sense 
    \begin{subequations} \label{eq:weak_snpp}
        \begin{align}
            \kappa \int_\Omega \partial_t \bv \cdot \bs{\varphi} \diff \X + \int_\Omega \nabla \bv : \nabla \bs{\varphi} \diff \X + \int_\Omega (c^+ - c^-) \nabla \psi \cdot \bs{\varphi} \diff \X &= 0, \label{eq:weak_Stokes}\\
            \langle \partial_t \cpm, \varphi \rangle - \int_\Omega \cpm \, \bv \cdot \nabla \varphi \diff \X + \int_\Omega \varMa \bigl( \nabla \cpm \pm \cpm \nabla \psi \bigr) \cdot \nabla \varphi \diff \X &= 0, \label{eq:weak_np}\\
            \int_\Omega \varMa \nabla \psi \cdot \nabla \varphi \diff \X + \int_\Gamma (\tau \psi - \xi) \varphi \dS - \int_\Omega (c^+ - c^-) \, \varphi \diff \X &= 0\label{eq:weak_Poisson}
        \end{align}
    \end{subequations}
    for a.e. $t \in (0,T)$, and for every $\bs{\varphi} \in W^{1,2}_{0,\sigma}(\Omega), \varphi \in W^{1,2}(\Omega)$.
\end{definition}

\begin{remark}\label{remark:p}
    Given a weak solution $(\bv,\cpm,\psi)$, then the pressure $p$ appearing in the strong formulation of the instationary Stokes equation \eqref{eq:s} can be reconstructed a posteriori up to an additive constant such that
    \begin{equation*}
          \nabla p  = - \kappa\partial_t \bv + \Delta \bv - (c^+ - c^-) \nabla \psi
    \end{equation*}
    holds in the sense of distributions (see, e.g., \cite[Chapter IV, Section 2.6]{sohr_01}). Since the right-hand side belongs to $L^2(0,T;L^2(\Omega))$, choosing $p$ such that $\int_\Omega p \diff \X = 0$, then $p \in L^2(0,T;W^{1,2}(\Omega))$.
\end{remark}
\noindent To proceed with the analysis of our optimal control problem, we start by proving the existence of weak solutions.
\begin{thm}[Existence of weak solutions]\label{thm:weak_existence}
    For every control $\xi \in \CXR$, there exists a global-in-time weak solution $(\bv, \cpm, \psi)$ of the SNPP system \eqref{eq:snpp}--\eqref{eq:ic} 
    in the sense of \zcref{definition:weak_sol}.
    Moreover, there exists a positive constant $C_{1,R}$ such that
    \begin{equation}\label{eq:main_norm_ineq}
        \begin{split}
            &\| \bv \|_{W^{1,2}(0,T;L^2(\Omega)) \cap L^2(0,T;W^{1,2}_{0,\sigma}(\Omega) \cap W^{2,2}(\Omega))}\\
            &\quad + \| \cpm \|_{W^{1,2}(0,T;(W^{1,2}(\Omega))^*) \cap L^\infty(0,T;L^2(\Omega)) \cap L^2(0,T;W^{1,2}(\Omega))} + \| \psi \|_{L^\infty(0,T;W^{1,6}(\Omega))} \leq C_{1,R}.
        \end{split}
    \end{equation}
    Finally, the charges are nonnegative and essentially bounded, i.e.,
    \begin{equation}\label{eq:bound_charges}
        0 \leq \cpm \leq M^{\pm}_R
    \end{equation}
    for a certain constant $M^{\pm}_R > 0$.
    In particular, $M^{\pm}_R$ is uniform with respect to all controls $\xi \in \CXR$.
\end{thm}
\noindent To obtain a well-defined optimal control problem, we need uniqueness of weak solutions. We establish it in the following theorem as a consequence of the Lipschitz continuous dependence of the solution on the initial data and the control.
\begin{thm}[Lipschitz continuous dependence and uniqueness] \label{thm:uniqueness}
    For every couple of controls $\xi_1, \xi_2 \in \CXR$ and initial data $(\bv_{0,1}, \cpm_{0,1}), (\bv_{0,2}$, $\cpm_{0,2})$ satisfying \ref{ass:initial_data} with associated weak solutions $(\bv_1, \cpm_1, \psi_1), (\bv_2, \cpm_2, \psi_2)$, 
    the following Lipschitz continuous dependence inequality holds
    \begin{equation}\label{eq:cont_dep_ineq}
         \begin{split}
            &\| \bv_1 - \bv_2 \|_{L^{\infty}(0,T;L^2(\Omega)) \cap L^2(0,T;W^{1,2}(\Omega))} + \sigmp{\| \cpm_1 - \cpm_2 \|_{L^{\infty}(0,T;L^2(\Omega)) \cap L^2(0,T;W^{1,2}(\Omega))} } \\
            &\quad\quad +  \| \psi_1 - \psi_2 \|_{ L^2(0,T;W^{1,2}(\Omega))}\\
            &\quad \leq C_{2,R} \Big( \| \bv_{0,1} - \bv_{0,2} \|_{L^2(\Omega)} +  \sigmp{\| \cpm_{0,1} - \cpm_{0,2} \|_{L^2(\Omega)}} + \| \xi_1 - \xi_2 \|_{L^2(0,T;L^2(\Gamma))}\Big),
        \end{split}
    \end{equation}
    for a positive constant $C_{2,R}$ (not depending on $\xi_1, \xi_2$). Moreover, there exists a constant $C_{3,R}>0$, such that
    \begin{equation}\label{eq:improved_cont_dep_ineq}
        \| \psi_1 - \psi_2 \|_{ L^{\infty}(0,T;W^{1,6}(\Omega))} \leq C_{3,R} \Big( \| \bv_{0,1} - \bv_{0,2} \|_{L^2(\Omega)} +  \sigmp{\| \cpm_{0,1} - \cpm_{0,2} \|_{L^2(\Omega)}} + \| \xi_1 - \xi_2 \|_{L^{\infty}(\Sigma)}\Big).
    \end{equation}
    In particular, the Cauchy problem \eqref{eq:snpp}--\eqref{eq:ic} admits a unique weak solution.
\end{thm}

\begin{remark}
    The constants $C_{1,R}, C_{2,R}, C_{3,R}, M^\pm_R, $ as highlighted by the subscript, depend on $R$, 
    but they also depend on the system parameters $\Omega, T, \bd, \tau, \kappa, \lambda = \varepsilon, M^\pm, \norm{\bv_0}_{W^{1,2}(\Omega)}$.
\end{remark}

\subsection{Existence of weak solutions (proof of \zcref{thm:weak_existence})}\label{sec:weak_exist}
We prove \zcref{thm:weak_existence} with a classical strategy. First, we introduce a regularized version of the system, depending on a parameter $\delta > 0$, for which the existence of solutions can be shown more easily. Then, we derive a priori estimates that are uniform with respect to $\delta$. Finally, these estimates give enough compactness to pass to the limit as $\delta \to 0$ and obtain a weak solution of \eqref{eq:snpp}--\eqref{eq:ic}.
The regularized system was introduced in \cite{Fischer_Saal_17} for a Navier--Stokes--Nernst--Planck--Poisson system with scalar diffusion coefficients and later used in \cite{homberg_lasarzik_plato_2025}, where the existence of global weak solutions to a similar anisotropic system is shown. While the regularization procedure is similar in the present case, the a priori estimates derived in \cite{homberg_lasarzik_plato_2025} are obtained under stronger regularity assumptions on the external potential $\xi$ and do not carry over to our setting, where only $\xi \in L^\infty(\Sigma)$ is assumed. Thus, in the following, we are going to only sketch the regularization part (see \zcref{subsubsect:reg_system}), omit the limit $\delta \to 0^+$, and we will focus on the a priori estimates, which are the main novelty (see \zcref{subsubsec:formal_apriori_est}). To keep the notation light, we derive them formally on the system \eqref{eq:snpp}--\eqref{eq:ic}.

\subsubsection{The regularized system}\label{subsubsect:reg_system}
 We briefly outline the regularization procedure of \cite{homberg_lasarzik_plato_2025, Fischer_Saal_17}, which carries over to the present setting with only minor modifications. In addition, we replace the director $\bd$ by a mollified version $\bd_\delta$ in the equation for the internal electric potential. The smooth coefficients this produces let us apply higher-order elliptic regularity results. The mollification parameter $\delta$ is sent to zero in the final limit passage. The regularization itself is carried out in three steps.\\
 
 \noindent \textit{1. Show the existence of weak solutions to the regularized Nernst--Planck--Poisson subsystem.}\\[0.1em]
 Let $\bv \in  L^\infty(0,T;W^{1,2}_{0,\sigma}(\Omega))$ be the given velocity field. We look at the corresponding Nernst--Planck--Poisson subsystem, which we regularize employing an elliptic regularization for the right-hand side of the internal Poisson equation $(P_i)$: the charge difference $c^+ - c^-$ is replaced by a smoothed version $\rho_\delta$, obtained by solving an additional elliptic problem with parameter $\delta > 0$. Since $\psie$ is given by $(P_e)$, we collect the transport due to the velocity field and to the external electric potential in the drift $\bs{b} \coloneqq \bv \mp \varMa \nabla \psie$, which is a known function belonging to $L^\infty(0,T;L^6(\Omega))$, by the Sobolev embedding $W^{1,2}(\Omega) \hookrightarrow L^6(\Omega)$ and \eqref{eq:est_psie}. For $\bd_\delta \in C^{2,1}(\overline{\Omega})$, we solve
 \begin{align*}
    \partial_t \cpm + \nabla \cdot (\cpm \bs{b}) - \nabla \cdot \left(\varMadelta (\nabla \cpm \pm \cpm \nabla \psii)\right) &= 0\\
    - \nabla \cdot (\varMadelta \nabla \psii) &= \rho_\delta\\
    \rho_\delta -  \delta \nabla \cdot (\varMadelta \nabla \rho_\delta) &= c^+ - c^-
\end{align*}
with boundary conditions
\begin{gather*}
    (-\cpm \bs{b} + \varMadelta (\nabla \cpm \pm \cpm \nabla \psii)) \cdot \bn = 0, \quad
    \varMadelta \nabla \psii \cdot \bn + \tau \psii = 0, \quad
    \varMadelta \nabla \rho_\delta \cdot \bn + \tau \rho_\delta = 0.
\end{gather*}
Heuristically, for $\delta \to 0^+$ we have $\rho_\delta \to  c^+ - c^-$ and we recover $(P_i)$. 
Moreover, for $\delta > 0$, $\rho_\delta$ is more regular than $c^+ - c^-$, and so is $\nabla \psii$, which makes the drift term $\cpm \nabla \psii$ easier to handle. By Schaefer's fixed point theorem \cite[Sec.~9, Thm.~4]{evans}, we then obtain the existence of a solution
\begin{align*}
    \cpm \in W^{1,2}(0,T;(W^{1,2}(\Omega))^*) \cap L^2(0,T;W^{1,2}(\Omega)), \quad \cpm \geq 0,
    \quad \psii \in C([0,T];W^{2,2}(\Omega)).
\end{align*}
The uniqueness follows from a standard test, 
see the proof of \cite[Prop.~4.3]{homberg_lasarzik_plato_2025} and \cite[Lem.~4.1]{Fischer_Saal_17} for more details.\\

\noindent \textit{2. Show the local-in-time existence of weak solutions to the fully coupled regularized system.}\\[0.1em]
 To couple the Nernst--Planck--Poisson subsystem with the Stokes equation, we define the fixed-point map
 \begin{equation*}
      \mathcal{H}_\delta: L^\infty(0,T;W^{1,2}_{0,\sigma}(\Omega)) \to L^\infty(0,T;W^{1,2}_{0,\sigma}(\Omega)),
 \end{equation*}
 mapping a velocity $\bv$ to
\[
    \mathcal{H}_\delta(\bv) = e^{-\frac{t}{\kappa} \mathcal{A}_2} \bv_0 - \int_0^t \tfrac{1}{\kappa} \; e^{- \frac{(t-s)}{\kappa} \mathcal{A}_2} \mathcal{P}_H((c^+ - c^-) \nabla \psi) \diff s,
\]
where $(\cpm, \psii)$ is the solution of Step~1 corresponding to $\bv$ and $\psi = \psii + \psie$. Here, $\mathcal{A}_2$ denotes the Stokes operator in $L^2(\Omega)$, $e^{-t \mathcal{A}_2}$ denotes the bounded analytic semigroup generated by the Stokes operator, and $\mathcal{P}_H$ denotes the Helmholtz projection in $L^2(\Omega)$ onto the divergence-free functions. In other words, $\mathcal{H}_\delta(\bv)$ is the solution of the instationary Stokes equation with the electric force computed from $\bv$. We can show that at least locally in time $\mathcal{H}_\delta$ is a contraction on some ball in $L^\infty(0,T;W^{1,2}_{0,\sigma}(\Omega))$, and Banach's fixed point theorem provides the existence of a fixed point of $\mathcal{H}_\delta$, which then is the solution to the fully coupled regularized system.\\

\noindent \textit{3. Show that these local solutions may be extended to arbitrary times $T>0$.}\\[0.1em]
As usual, the extension argument relies on the a priori estimates showing that the solutions do not blow up in finite time.

\subsubsection{Formal a priori estimates} \label{subsubsec:formal_apriori_est}
The a priori estimates necessary for the existence of weak solutions can be proved in five steps:
\begin{enumerate}[label =(\arabic*)]
    \item \label{enum:nonnegativity_charges}
    \textbf{Nonnegativity of the charges:} $\cpm \geq 0$;
    \item \label{enum:energy_est}
    \textbf{Energy estimate:}
        $\bv \in L^\infty(0,T;L^2(\Omega)) \cap L^2(0,T;W^{1,2}(\Omega))$, 
        $\cpm \in L^\infty(0,T;L^1(\Omega))$,
        $\ssqrt{\cpm} \in L^2(0,T;W^{1,2}(\Omega))$, and
        $\psi \in L^\infty(0,T;W^{1,2}(\Omega)) \cap L^2(0,T;W^{1,6}(\Omega))$;
    \item \label{enum:grad_est_charges}
    \textbf{$L^2$-gradient estimate for the charges:}
        $\cpm \in L^\infty(0,T;L^2(\Omega)) \cap L^2(0,T;W^{1,2}(\Omega))$,
        which yields $\psi \in L^\infty(0,T;W^{1,6}(\Omega))$ by elliptic regularity;
    \item \label{enum:velocity_reg}
    \textbf{Regularity estimate for the velocity:}
        The electric volume force satisfies $(c^+ - c^-) \nabla \psi \in L^2(0,T;L^2(\Omega))$ by the previous step,
        hence $\bv \in W^{1,2}(0,T;L^2(\Omega)) \cap L^2(0,T;W^{2,2}(\Omega)) \hookrightarrow L^2(0,T;L^\infty(\Omega))$.
    \item \label{enum:time_reg_charges}
    \textbf{Time regularity estimate for the charges:}
    $\partial_t \cpm \in L^2(0,T;(W^{1,2}(\Omega))^*)$.
\end{enumerate}

\noindent Finally, to conclude the proof of \zcref{thm:weak_existence}, we will prove the boundedness of the charges to get \eqref{eq:bound_charges}. This will be done in the following final step:

\begin{enumerate}[resume, label =(\arabic*)]
    \item \textbf{Boundedness of the charges:}
    $\cpm \in L^{\infty}(Q)$.         \label{enum:bound_charges}
\end{enumerate}
Notice that the existence of weak solutions follows from the steps \ref{enum:nonnegativity_charges}--\ref{enum:time_reg_charges}. Thus, it is proved without a uniform $ L^\infty$-bound on $\cpm$, which will be derived only at the end.
All constants depend only on the problem data and on $R$ (unless stated and discussed otherwise).

\paragraph{Step \ref{enum:nonnegativity_charges}: Nonnegativity of the charges.}
After introducing the notation  $(\cpm)_- \coloneqq \max\{0, - \cpm \}$ to denote the negative part of $\cpm$, we test the Nernst--Planck equation \eqref{eq:weak_np} with $-(\cpm)_-$ (which is an admissible test function). Since $\bv$ vanishes at the boundary and $\nabla \cdot \bv = 0$, integrating by parts in the velocity term leads to
\begin{equation*}
    - \int_\Omega (\cpm)_- \, \bv \cdot \nabla (\cpm)_- \diff \X = - \frac{1}{2}\int_\Omega  \bv \cdot \nabla [(\cpm)_-]^2 \diff \X = \frac{1}{2}\int_\Omega \nabla \cdot \bv \, [(\cpm)_-]^2 \diff \X = 0.
\end{equation*}
We obtain
\begin{equation}\label{eq:noneg_0}
    \begin{split}
        &\frac{1}{2} \frac{\diff}{\diff t} \int_\Omega |(\cpm)_-|^2 \diff \X + \int_\Omega |\nabla (\cpm)_-|_{\varMa}^2 \diff \X\\
        &\quad = \mp  \int_\Omega (\cpm)_- \varMa \nabla \psi \cdot \nabla (\cpm)_- \diff \X \leq
        C \norm{(\cpm)_-}_{L^3(\Omega)} \norm{\nabla \psi}_{L^6(\Omega)} \norm{\nabla(\cpm)_-}_{L^2(\Omega)},
    \end{split}
\end{equation}
where we have exploited the boundedness of $\varMa$ given by assumption \ref{ass:anisotropy} and the H\"older inequality to estimate the right-hand side. Then, we use the Gagliardo--Nirenberg inequality \eqref{eq:GN_L3} to handle the $L^3$-norm, and finally the Young inequality with a small constant $\delta>0$ (fixed but yet to be chosen). We get:
\begin{align*}
 &\frac{1}{2} \frac{\diff}{\diff t} \int_\Omega |(\cpm)_-|^2 \diff \X + \int_\Omega |\nabla (\cpm)_-|_{\varMa}^2 \diff \X\\
 &\quad \leq C \Big(\norm{(\cpm)_-}^{1/2}_{L^2(\Omega)}\norm{\nabla(\cpm)_-}^{1/2}_{L^2(\Omega)} + \norm{(\cpm)_-}_{L^2(\Omega)}\Big) \norm{\nabla \psi}_{L^6(\Omega)} \norm{\nabla(\cpm)_-}_{L^2(\Omega)}\\
 &\quad \leq  C \norm{(\cpm)_-}^{1/2}_{L^2(\Omega)} \norm{\nabla \psi}_{L^6(\Omega)} \norm{\nabla(\cpm)_-}^{3/2}_{L^2(\Omega)} + C\norm{(\cpm)_-}_{L^2(\Omega)}\norm{\nabla \psi}_{L^6(\Omega)} \norm{\nabla(\cpm)_-}_{L^2(\Omega)}\\
 &\quad \leq \delta \norm{\nabla(\cpm)_-}_{L^2(\Omega)}^2 + C_{\delta} \Big(\norm{\nabla \psi}^4_{L^6(\Omega)} + \norm{\nabla \psi}^2_{L^6(\Omega)} \Big)\norm{(\cpm)_-}_{L^2(\Omega)}^2.
\end{align*}
Employing $|\cdot|_{\varMa} \geq |\cdot|$, choosing $\delta = 1/2$, and integrating in time over $(0,s)$ for $s \in (0,T)$ we end up with: 
\begin{equation}\label{eq:noneg_3}
    \begin{split}
         &\int_\Omega |(\cpm)_-|^2 \diff \X + \int_0^s\int_\Omega |\nabla (\cpm)_-|^2 \diff \X  \diff t\\
         &\quad \leq \int_\Omega |(\cpm_0)_-|^2 \diff \X  + C \int_0^s\Big(\norm{\nabla \psi}^4_{L^6(\Omega)} + \norm{\nabla \psi}^2_{L^6(\Omega)} \Big) \|(\cpm)_-\|_{L^2(\Omega)}^2\diff t.
    \end{split}
\end{equation}
By the Gronwall inequality, we conclude that
\begin{equation*}
    \int_\Omega |(\cpm)_-|^2 \diff \X + \int_0^s\int_\Omega |\nabla (\cpm)_-|_{\varMa}^2 \diff \X  \diff t \leq \bigg( \int_\Omega |(\cpm_0)_-|^2 \diff \X \bigg)e^{C \int_0^s (\norm{\nabla \psi}^4_{L^6(\Omega)} + \norm{\nabla \psi}^2_{L^6(\Omega)}) \diff t} = 0
\end{equation*}
where the last equality holds because, by hypothesis \ref{ass:initial_data}, $\cpm_0$ is nonnegative. Thus, $\cpm$ is nonnegative, and nonnegativity is preserved under the passage to the limit via weak $L^2(0,T;L^2(\Omega))$ convergence.
Before moving on, a final comment. 
The constant $C$ appearing in \eqref{eq:noneg_3} depends on $\|\cpm\|_{L^{\infty}(0,T;L^2(\Omega))}$ and $\|\psi\|_{L^{\infty}(0,T;W^{1,2}(\Omega))}$, which is atypical when deriving uniform a priori estimates. However, in this case, this is not an issue, since in the end the term where it appears vanishes. The only important aspect to perform this estimate at the approximate level is that the approximate solution has (at least) this regularity, but this is natural and anyway required, since it is the regularity we recover in the limit. 
 
\paragraph{Step \ref{enum:energy_est}: Energy estimate.} 
Testing the Nernst--Planck equation \eqref{eq:weak_np} with $\varphi = 1$ and using the no-flux boundary condition gives mass conservation
\begin{equation}\label{eq:const_mean}
    \frac{\diff}{\diff t} \int_\Omega \cpm \diff \X = 0,
    \qquad\text{hence}\qquad
    \int_\Omega \cpm(t) \diff \X = \int_\Omega \cpm_0 \diff \X.
\end{equation}
Together with $\cpm \geq 0$ this yields that $\cpm$ are uniformly bounded in $ L^\infty(0,T;L^1(\Omega))$.
The energy estimate is obtained by the testing scheme
\(
    S \cdot \bv
    + NP \cdot (\ln \cpm + 2 \pm \psii)
    + \partial_t (P_i) \cdot \psii,
\)
where we deliberately test the Nernst--Planck equation with the internal potential $\psii$ only. Here and in the following we are using the splitting $\psi = \psii + \psie$ introduced in \eqref{eq:splitting}, and we will employ the estimates \eqref{eq:est_psii} and \eqref{eq:est_psie}. 
The goal is the energy-dissipation inequality
\begin{equation}
    \begin{split}\label{eq:energy_estimate}
        &\left[\int_\Omega \Big(\frac{\kappa}{2} |\bv|^2 + \sigmp{\cpm (\ln \cpm + 1)} + \frac{1}{2}|\nabla \psii|_{\varMa}^2 \Big) \diff \X + \frac{\tau}{2} \int_\Gamma |\psii|^2 \dS \right]_0^t\\
        &\quad + \int_0^t \int_\Omega \Big(|\nabla \bv|^2 + 4 \sigmp{|\nabla \ssqrt{\cpm}|_{\varMa}^2} + (c^+ + c^-) |\nabla \psii|_{\varMa}^2 + 2|c^+ - c^-|^2 \Big)\diff \X \diff t \leq C.
    \end{split}
\end{equation}
Notice that here $\psi_i(0)$ is not a prescribed initial datum, but is the solution of $(P_i)$ having as right-hand side $c^+(0) - c^-(0) = c^+_0 - c^-_0 \in L^{\infty}(\Omega) \subseteq L^2(\Omega)$. 
Thus, by \eqref{eq:est_psii}, both the gradient and the boundary term of $\psii$ at time zero in \eqref{eq:energy_estimate} can be controlled by the $L^2(\Omega)$-norm of $\cpm_0$.
Testing the Stokes equation \eqref{eq:weak_Stokes} with $\bv$ gives
\begin{equation}\label{eq:energy_NS}
    \frac{\diff}{\diff t} \frac{\kappa}{2} \int_\Omega |\bv|^2 \diff \X
    + \int_\Omega |\nabla \bv|^2 \diff \X
    = - \int_\Omega (c^+ - c^-) \nabla \psi \cdot \bv \diff \X = - \int_\Omega (c^+ - c^-) (\nabla \psii + \nabla \psie) \cdot \bv \diff \X.
\end{equation}
Then, we test the Nernst--Planck equation \eqref{eq:weak_np} with $\ln \cpm + 2 \pm \psii$ and sum over the $\pm$ indices. 
Notice that, since $\cpm \geq 0$ is only nonnegative and not bounded away from $0$, $\ln{\cpm}$ is not an admissible test function. 
This issue can be rigorously addressed by testing with $\ln{(\cpm + \delta)}$ for a constant $\delta > 0$ and then passing to the limit as $\delta \to 0^+$ (see \cite[pp.~16--17]{homberg_lasarzik_plato_2025}). In the following, we proceed formally, omitting the parameter $\delta$. We obtain:
\begin{align}
    &\frac{\diff}{\diff t} \int_\Omega \sigmp{\cpm (\ln \cpm + 1)} \diff \X + \int_\Omega \partial_t (c^+ - c^-) \psii \diff \X - \int_\Omega (c^+ - c^-) \bv \cdot \nabla \psii \diff \X \nonumber\\
    &\qquad + \int_\Omega \sigmp{\varMa \bigl( \nabla \cpm \pm \cpm \nabla \psi \bigr) \cdot \bigl( \nabla \ln \cpm \pm \nabla \psii \bigr)} \diff \X
    = 0, \label{eq:energy_NP}
\end{align}
where the convective contribution $\int_\Omega \cpm \bv \cdot \nabla \ln \cpm \diff \X$ vanishes by $\nabla \cdot \bv = 0$ and $\bv|_{\Gamma} = 0$.
Expanding the dissipation integrand and using $\psi = \psii + \psie$ (together with $\lamMa = \varMa$ from hypothesis \ref{ass:struct_cond}, see the \zcref{remark:lam_var}), we obtain
\begin{align*}
    &\sigmp{\varMa \bigl( \nabla \cpm \pm \cpm \nabla \psi \bigr) \cdot \bigl( \nabla \ln \cpm \pm \nabla \psii \bigr)}\\
    &\quad  = \sigmp{\varMa \frac{\nabla \cpm}{\ssqrt{\cpm}} \cdot \frac{\nabla \cpm}{\ssqrt{\cpm}} \pm \varMa \nabla \cpm \cdot \nabla \psii} + \sigmp{\pm \cpm \,\varMa \nabla \psie \cdot \bigl( \nabla \ln \cpm \pm \nabla \psii \bigr)}\\
    &\quad \quad+ \sigmp{\pm \cpm \,\varMa \nabla \psii \cdot \bigl( \nabla \ln \cpm \pm \nabla \psii \bigr)}\\
    &\quad
    = 4 \sigmp{|\nabla \ssqrt{\cpm}|_{\varMa}^2}
    + (c^+ + c^-) |\nabla \psii|_{\varMa}^2
    + \sigmp{\pm \cpm \,\varMa \nabla \psie \cdot \bigl( \nabla \ln \cpm \pm \nabla \psii \bigr)}\\
    &\quad \quad 
    + \underbrace{\sigmp{\pm \varMa \nabla \cpm \cdot \nabla \psii \pm \cpm \,\varMa \nabla \psii \cdot \nabla \ln \cpm}}_{=: I_1}.
\end{align*}
The cross term $I_1$, integrated on $\Omega$, is rewritten via the internal Poisson equation $(P_i)$:
\begin{align*}
    \int_{\Omega }I_1 \diff \X
    &= 2 \int_\Omega \varMa \nabla \psii \cdot \nabla (c^+ - c^-) \diff \X
    = 2 \int_\Omega |c^+ - c^-|^2 \diff \X
    - 2 \int_\Gamma \tau \psii (c^+ - c^-) \dS.
\end{align*}
Differentiating $(P_i)$ in time and testing with $\psii$ gives
\begin{align}
    \int_\Omega \partial_t (c^+ - c^-) \psii \diff \X
    &= \int_\Omega \varMa \nabla \partial_t \psii \cdot \nabla \psii \diff \X
    + \tau \int_\Gamma \partial_t \psii \, \psii \dS \nonumber\\
    &= \frac{1}{2} \frac{\diff}{\diff t}
        \left( \int_\Omega |\nabla \psii|_{\varMa}^2 \diff \X + \tau \int_\Gamma |\psii|^2 \dS \right), \label{eq:energy_Pi}
\end{align}
where we exploited the fact that the director $\bd$ does not depend on time.
Adding \eqref{eq:energy_NS}--\eqref{eq:energy_Pi}, the transport term $\int_\Omega (c^+ - c^-) \nabla \psii \cdot \bv \diff \X$ cancels against the corresponding contribution from Stokes, and we arrive at
\begin{align}
    &\frac{\diff}{\diff t} \Biggl(
        \frac{\kappa}{2} \int_\Omega |\bv|^2 \diff \X
        + \int_\Omega \sigmp{\cpm (\ln \cpm + 1)} \diff \X
        + \frac{1}{2} \int_\Omega |\nabla \psii|_{\varMa}^2 \diff \X
        + \frac{\tau}{2} \int_\Gamma |\psii|^2 \dS
    \Biggr) \nonumber\\
    &\qquad
    + \int_\Omega \Big( |\nabla \bv|^2
        + 4 \sigmp{|\nabla \ssqrt{\cpm}|_{\varMa}^2}
        + (c^+ + c^-) |\nabla \psii|_{\varMa}^2
        + 2 |c^+ - c^-|^2 \Big)\diff \X
    = I_2 + I_3 + I_4, \label{eq:energy_balance}
\end{align}
where the right-hand side consists of
\begin{align*}
    I_2
    &\coloneqq - \int_\Omega (c^+ - c^-) \bv \cdot \nabla \psie \diff \X,\\
    I_3
    &\coloneqq - \int_\Omega \sigmp{\pm \cpm \,\varMa \nabla \psie \cdot \bigl( \nabla \ln \cpm \pm \nabla \psii \bigr)} \diff \X,\\
    I_4
    &\coloneqq 2 \int_\Gamma \tau \psii (c^+ - c^-) \dS.
\end{align*}

\smallskip
\noindent\textit{Estimate of $|I_2|$.}
By H\"older's inequality,
\[
    |I_2|
    \leq \norm{c^+ - c^-}_{L^2(\Omega)} \norm{\bv}_{L^3(\Omega)} \norm{\nabla \psie}_{L^6(\Omega)}.
\]
The Young inequality and the Gagliardo--Nirenberg inequality \eqref{eq:GN_L3_0} yield, for any $\delta > 0$ and some constant $C > 0$ depending on $\delta$:
\begin{align*}
    |I_2|
    &\leq \delta \norm{c^+ - c^-}_{L^2(\Omega)}^2
        + C \norm{\bv}_{L^3(\Omega)}^2 \norm{\nabla \psie}_{L^6(\Omega)}^2\\
    &\leq \delta \norm{c^+ - c^-}_{L^2(\Omega)}^2
        + C \norm{\bv}_{L^2(\Omega)} \norm{\nabla \bv}_{L^2(\Omega)} \norm{\nabla \psie}_{L^6(\Omega)}^2\\
    &\leq \delta \norm{c^+ - c^-}_{L^2(\Omega)}^2 
        + \delta \norm{\nabla \bv}_{L^2(\Omega)}^2
        + C \norm{\bv}_{L^2(\Omega)}^2,
\end{align*}
where we have also used the estimate \eqref{eq:est_psie} according to which $\nabla \psie$ is uniformly bounded in $L^\infty(0,T;L^6(\Omega))$. The first two terms on the right are absorbed into the left-hand side of \eqref{eq:energy_balance}.

\smallskip
\noindent\textit{Estimate of $|I_3|$.}
Rewriting $\nabla \ln \cpm = 2 \nabla \ssqrt{\cpm} / \ssqrt{\cpm}$ we obtain
\begin{align*}
    |I_3|
    &= \Biggl|
        \int_\Omega \Big( \sigmp{\pm 2 \ssqrt{\cpm} \,\varMa \nabla \psie \cdot \nabla \ssqrt{\cpm}}
        + (c^+ + c^-) \varMa \nabla \psie \cdot \nabla \psii \Big) \diff \X
    \Biggr|\\
    &\leq
        C \sigmp{\norm{\ssqrt{\cpm}}_{L^3(\Omega)} \norm{\nabla \psie}_{L^6(\Omega)} \norm{\nabla \ssqrt{\cpm}}_{L^2(\Omega)}}\\
    &\quad + C \norm{\nabla \psie}_{L^{\infty}(0,T;L^6(\Omega))} \sigmp{\norm{\ssqrt{\cpm} \nabla \psii}_{L^2(\Omega)} \norm{\ssqrt{\cpm}}_{L^3(\Omega)}}.
\end{align*}
By the uniform estimate \eqref{eq:est_psie}, the Gagliardo--Nirenberg inequality \eqref{eq:GN_L3} applied for $f = \ssqrt{\cpm}$, and Young's inequality, we obtain
\begin{align*}
    |I_3|
    &\leq \delta \sigmp{\norm{\nabla \ssqrt{\cpm}}_{L^2(\Omega)}^2} + \delta \sigmp{\norm{\ssqrt{\cpm} \nabla \psii}_{L^2(\Omega)}^2} + C \sigmp{\norm{\ssqrt{\cpm}}_{L^2(\Omega)}^2}.
\end{align*}
The first two contributions are absorbed into the dissipation on the left-hand side of \eqref{eq:energy_balance}.
For the Gronwall term, we use the elementary bound
\[
    \norm{\ssqrt{\cpm}}_{L^2(\Omega)}^2
    = \int_\Omega \cpm \diff \X = \int_\Omega \cpm_0 \diff \X.
\]

\smallskip
\noindent\textit{Estimate of $|I_4|$.}
We use the embedding
\begin{equation}\label{eq:Linf_W14}
    \norm{f}_{L^\infty(\Omega)}
    \leq C \norm{f}_{W^{1,4}(\Omega)}
    \leq C \Big(\norm{\nabla f}_{L^2(\Omega)}^{1/4} \norm{\nabla f}_{W^{1,2}(\Omega)}^{3/4} + \norm{f}_{W^{1,2}(\Omega)} \Big),
\end{equation}
obtained by the continuous embedding $W^{1,4}(\Omega) \hookrightarrow L^{\infty}(\Omega)$ and the Gagliardo--Nirenberg interpolation inequality \eqref{eq:GN_L4}. Additionally, we use the trace inequality \eqref{eq:trace_L2}.
Then,
\begin{align*}
    |I_4|
    &\leq C \norm{\psii}_{L^\infty(\Omega)} \sigmp{\norm{\ssqrt{\cpm}}_{L^2(\Gamma)}^2}\\
    &\leq C \Biggl(
        \norm{\psii}_{W^{1,2}(\Omega)}
        + \norm{\nabla \psii}_{L^2(\Omega)}^{1/4} \norm{\psii}_{W^{2,2}(\Omega)}^{3/4}
    \Biggr)
    \sigmp{\norm{\ssqrt{\cpm}}_{L^2(\Omega)} \norm{\ssqrt{\cpm}}_{W^{1,2}(\Omega)}}\\
    &\leq C \Biggl(
        \norm{\psii}_{W^{1,2}(\Omega)}
        + \norm{\nabla \psii}_{L^2(\Omega)}^{1/4} \norm{c^+ - c^-}_{L^2(\Omega)}^{3/4}
    \Biggr)
    \sigmp{ 1 + \norm{\nabla \ssqrt{\cpm}}_{L^2(\Omega)}},
\end{align*}
where we used that $\snorm{\ssqrt{\cpm}}_{L^2(\Omega)}$ is bounded by the $L^\infty(0,T;L^1)$-bound on $\cpm$ and the elliptic regularity estimates \eqref{eq:est_psii}.
Applying Young's inequality first with conjugate exponents $(p,p') = (2, 2)$ and then with $(4, 4/3)$, we conclude for any $\delta > 0$
\begin{align*}
    |I_4|
    &\leq
        C \Biggl(\norm{\psii}_{W^{1,2}(\Omega)}^2 + \norm{\nabla \psii}_{L^2(\Omega)}^{1/2} \norm{c^+ - c^-}_{L^2(\Omega)}^{3/2} \Biggr) 
        + \delta \Biggl( 1 + \sigmp{\norm{\nabla \ssqrt{\cpm}}_{L^2(\Omega)}^2}\Biggr)\\
    &\leq
        C \norm{\psii}_{W^{1,2}(\Omega)}^2
        + \delta \norm{c^+ - c^-}_{L^2(\Omega)}^2
        + \delta \sigmp{\norm{\nabla \ssqrt{\cpm}}_{L^2(\Omega)}^2} + C.
\end{align*}
The gradient and charge-difference terms are absorbed into the left-hand side of \eqref{eq:energy_balance},
while $\norm{\psii}_{W^{1,2}(\Omega)}^2$ can be bounded using the generalized Poincaré inequality \eqref{eq:gen_Poincare_ineq} and then handled with Gronwall.

Collecting the estimates for $I_2$--$I_4$, choosing $\delta > 0$ sufficiently small, and applying Gronwall's inequality \cite[Appx.~B.2.j]{evans}, we obtain \eqref{eq:energy_estimate}.
In particular, since the function $s (\ln{s} + 1)$ for $s\in(0,+\infty)$ is bounded from below, we get
\[
    \bv \in L^\infty(0,T;L^2(\Omega)) \cap L^2(0,T;W^{1,2}(\Omega)),
    \quad
    \ssqrt{\cpm} \in L^2(0,T;W^{1,2}(\Omega)),
    \quad
    \psii \in L^\infty(0,T;W^{1,2}(\Omega)).
\]
Moreover, since we have also shown that $c^+ - c^-$ is uniformly bounded in $L^2(0,T;L^2(\Omega))$,  we find that $\psii$ is uniformly bounded in $L^2(0,T;W^{2,2}(\Omega)) \hookrightarrow L^2(0,T;W^{1,6}(\Omega))$ by elliptic regularity estimate \eqref{eq:est_psii},
and by estimate \eqref{eq:est_psie} we find $\psie$ in $L^{\infty}(0,T;W^{1,6}(\Omega))$. Thus, we conclude
\begin{equation*}
    \| \psi \|_{L^2(0,T;W^{1,6}(\Omega))} \leq \| \psii \|_{L^2(0,T;W^{1,6}(\Omega))} + \| \psie \|_{L^2(0,T;W^{1,6}(\Omega))} \leq C.
\end{equation*}
In particular, $\psi$ is also uniformly bounded in  $L^2(0,T;W^{1,6}(\Omega)) \hookrightarrow L^2(0,T;L^\infty(\Omega))$.

\paragraph{Step \ref{enum:grad_est_charges}: $L^2$-gradient estimate for the charges.}
Testing the Nernst--Planck equation \eqref{eq:weak_np} with $\cpm$ and using $\nabla \cdot \bv = 0$ together with the no-flux boundary condition, we obtain
\begin{align*}
    &\frac{1}{2} \frac{\diff}{\diff t} \int_\Omega |\cpm|^2 \diff \X
    + \int_\Omega \Big(|\nabla \cpm|_{\varMa}^2 
    \pm \frac{1}{2} \varMa \nabla \psi \cdot \nabla(\cpm)^2 \Big)\diff \X\\
    &\quad = 
    \frac{1}{2} \frac{\diff}{\diff t} \int_\Omega |\cpm|^2 \diff \X
    + \int_\Omega \Big(|\nabla \cpm|_{\varMa}^2 
    \pm \frac{1}{2} (c^+ - c^-) (\cpm)^2 \Big)\diff \X \pm \frac{1}{2} \int_\Gamma (\xi - \tau \psi) (\cpm)^2\dS
    = 0,
\end{align*}
where we used the Poisson equation \eqref{eq:weak_Poisson} tested with $(\cpm)^2$ (thanks to the Assumption \ref{ass:struct_cond}) to rewrite the last term.
 Notice that $(\cpm)^2$ is not a priori an admissible test function since we have not proved yet that $\cpm$ is bounded. Thus, a rigorous approach would require introducing a truncation of $\cpm$
    \begin{equation*}
        \cpm_n \coloneqq \begin{cases}
            \cpm &\quad  0 \leq \cpm < n,\\
            n & \quad \cpm  \geq n,
        \end{cases}
    \end{equation*}
doing all the testing with it, and, in the end, passing to the limit as $n \to +\infty$. We will not do it since it is standard and we do not want to overload the notation.
We now add the equations for $c^+$ and $c^-$ and use the fact that
\[
    \int_\Omega (c^+ - c^-) ((c^+)^2 - (c^-)^2) \diff \X \geq 0
\]
which follows from the monotonicity of the map $x \mapsto x^2$ to get
\begin{align*}
    &\frac{1}{2} \frac{\diff}{\diff t} \int_\Omega \sigmp{|\cpm|^2} \diff \X
    + \int_\Omega \sigmp{|\nabla \cpm|_{\varMa}^2 } \diff \X 
    \leq
    C \left( \norm{\xi}_{L^{\infty}(\Sigma)} + \norm{\psi}_{L^\infty(\Omega)} \right) \sigmp{\norm{\cpm}_{L^2(\Gamma)}^2}\\
    &\quad \leq 
        \frac{1}{2} \int_\Omega \sigmp{|\nabla \cpm|^2} \diff \X 
    + C \left( 1 + \norm{\xi}^2_{L^{\infty}(\Sigma)} + \norm{\psi}^2_{L^\infty(\Omega)} \right) \sigmp{\norm{\cpm}_{L^2(\Omega)}^2}.
\end{align*}
where we used the trace estimate \eqref{eq:trace_L2}.
Absorbing the gradient term and applying Gronwall's inequality, we conclude that 
\begin{equation}\label{eq:grad_est}
    \| c^{\pm} \|_{L^\infty(0,T;L^2(\Omega))} + \| \nabla c^{\pm} \|_{L^2(0,T;L^2(\Omega))} \leq C
\end{equation}
since by the energy estimate from Step \ref{enum:energy_est} we have $\psi\in L^2(0,T;L^\infty(\Omega))$.
Elliptic regularity \eqref{eq:est_psii} for the internal Poisson equation $(P_i)$ then upgrades the potential estimate to
\begin{equation} \label{eq:psi_est}
     \| \psi \|_{L^\infty(0,T;W^{1,6}(\Omega))} \leq C.
\end{equation}

\paragraph{Step \ref{enum:velocity_reg}: Regularity estimate for the velocity.}
From Step~\ref{enum:grad_est_charges} we have $\cpm \in L^2(0,T;L^6(\Omega))$ by the Sobolev embedding $W^{1,2}(\Omega) \hookrightarrow L^6(\Omega)$,  and
$\nabla \psi \in L^\infty(0,T;L^6(\Omega))$. Therefore, the electric volume force satisfies
\[
    (c^+ - c^-) \nabla \psi
    \in L^2\bigl(0,T; L^6(\Omega)\bigr) \cdot L^\infty\bigl(0,T; L^6(\Omega)\bigr) \hookrightarrow L^2(0,T;L^3(\Omega))
    \hookrightarrow L^2(0,T;L^2(\Omega)).
\]
Maximal $L^2$-regularity for the nonstationary Stokes system \cite[Section IV, Thm.~2.7.1]{sohr_01} yields
\[
    \bv \in W^{1,2}(0,T;L^2(\Omega)) \cap L^2(0,T;W^{2,2}(\Omega))
    \hookrightarrow L^2(0,T;L^\infty(\Omega)).
\]

\paragraph{Step \ref{enum:time_reg_charges}: Time regularity estimate for the charges.}
We rewrite the Nernst--Planck equation \eqref{eq:np} as
\[
    \partial_t \cpm
    = - \nabla \cdot (\cpm \bv)
    - \nabla \cdot \bigl( \varMa (\nabla \cpm \pm \cpm \nabla \psi) \bigr).
\]
The fluxes belong to $L^2(0,T;L^2(\Omega))$: indeed,
$\nabla \cpm \in L^2(0,T;L^2(\Omega))$ by Step~\ref{enum:grad_est_charges},
$\cpm \nabla \psi \in L^2(0,T;L^2(\Omega))$ by the same product argument as in Step \ref{enum:velocity_reg}, and
$\cpm \bv \in L^2(0,T;L^2(\Omega))$ since $\cpm \in L^\infty(0,T;L^2(\Omega))$ and
$\bv \in L^2(0,T;L^\infty(\Omega))$ by Step~\ref{enum:velocity_reg}.
Consequently,
\[
    \partial_t \cpm \in L^2(0,T;(W^{1,2}(\Omega))^*),
\]
which completes the regularity required in \zcref{definition:weak_sol}.

\paragraph{Step \ref{enum:bound_charges}: Boundedness of the charges.} 
The argument follows a Moser--Alikakos iteration \cite[Thm.~3.1]{alikakos_79} along the exponents $p_k \coloneqq 2^k$ for $k \geq 1$, so that $p_1 = 2$ and $p_k = 2 p_{k-1}$ for $k \geq 2$.
For $k \geq 1$, we abbreviate
\[
    A_k  = \norm{\cpm}_{L^\infty(0,T;L^{p_k}(\Omega))} ,
\]
and we write $C$ for constants that may change from line to line but never depend on $k$.
We recursively prove the following uniform bound on $(A_k)_{k \geq 2}$ 
\begin{align}\label{eq:ih}
    A_k \leq \max\left(A_1, \norm{\cpm_0}_{L^\infty(\Omega)} \right) \prod_{j=2}^k \bigl( C \, (p_j^{10} + 1) \bigr)^{1/p_j},
\end{align}
and finally pass to the limit $k \to +\infty$. Notice that $A_1 = \norm{\cpm}_{L^\infty(0,T;L^2(\Omega))}$ and in \eqref{eq:grad_est} we have already shown that this quantity is uniformly bounded by a constant depending only on the problem's data and on $R$.
We test the Nernst--Planck equation \eqref{eq:weak_np} with $p_k (\cpm)^{p_k - 1}$.  
Notice that we again do not know that this is an admissible test function, and a rigorous approach would require truncating $c^{\pm}$ and then, in the end, passing to the limit in the truncation. Since this is standard, to avoid overloading the notation and an additional limiting procedure, we will proceed formally. 
The convective term vanishes,
\[
    p_k \int_\Omega \cpm \, \bv \cdot \nabla (\cpm)^{p_k - 1} \diff \X
    = (p_k - 1) \int_\Omega \bv \cdot \nabla (\cpm)^{p_k} \diff \X = 0,
\]
by $\nabla \cdot \bv = 0$ and the no-slip boundary condition.
Since 
\begin{align*}
    p_k \varMa \nabla \cpm \cdot \nabla ((\cpm)^{p_k - 1}) &= p_k (p_k - 1) (\cpm)^{p_k-2} \varMa \nabla \cpm \cdot \nabla \cpm 
    = p_k (p_k - 1) \left| (\cpm)^{p_k/2 - 1} \nabla \cpm \right|_{\varMa}^2\\
    &= p_k (p_k - 1) \frac{4}{p_k^2} |\nabla (\cpm)^{p_k/2}|_{\varMa}^2
    =  (p_k - 1) \frac{4}{p_k} |\nabla (\cpm)^{p_k/2}|_{\varMa}^2
\end{align*}
and 
\begin{align*}
    p_k \cpm \nabla (\cpm)^{p_k - 1} &= p_k (p_k -1) (\cpm)^{p_k - 1} \nabla \cpm = p_k (p_k - 1) (\cpm)^{p_k/2} \frac{2}{p_k} \nabla (\cpm)^{p_k/2}\\
    &=2 (p_k -1) (\cpm)^{p_k/2} \nabla (\cpm)^{p_k/2},
\end{align*}
we obtain
\begin{align}\label{eq:moser_iteration_1}
    \hspace{-1em}\frac{\diff}{\diff t} \int_\Omega (\cpm)^{p_k} \diff \X
    + \frac{4 (p_k-1)}{p_k} \int_\Omega |\nabla (\cpm)^{p_k/2}|_{\varMa}^2 \diff \X
    = \mp \, 2 (p_k - 1) \int_\Omega (\cpm)^{p_k/2} \, \varMa \nabla \psi \cdot \nabla (\cpm)^{p_k/2} \diff \X .
\end{align}
Then, we employ the H\"older inequality on the right-hand side, obtaining
\begin{align*}
    2 (p_k - 1) \left| \int_\Omega (\cpm)^{p_k/2} \, \varMa \nabla \psi \cdot \nabla (\cpm)^{p_k/2} \diff \X \right|
    &\leq C p_k \norm{\nabla \psi}_{L^6(\Omega)} \norm{(\cpm)^{p_k/2}}_{L^3(\Omega)} \norm{\nabla (\cpm)^{p_k/2}}_{L^2(\Omega)}\\
    &\leq C p_k \norm{(\cpm)^{p_k/2}}_{L^3(\Omega)} \norm{\nabla (\cpm)^{p_k/2}}_{L^2(\Omega)},
\end{align*}
where we also have used the fact that $\psi$ is bounded in $L^\infty(0,T;W^{1,6}(\Omega))$ by \eqref{eq:psi_est}.
Using the Gagliardo--Nirenberg interpolation inequality \eqref{eq:gagliardo_nirenberg_L1} from \zcref{lem:gagliardo_nirenberg_L1} to handle the $L^3$-norm, we get
 \begin{align*}
    &2 (p_k - 1) \left| \int_\Omega (\cpm)^{p_k/2} \, \varMa \nabla \psi \cdot \nabla (\cpm)^{p_k/2} \diff \X \right|\\
    &\quad \leq C p_k \Big(\norm{\nabla (\cpm)^{p_k/2}}_{L^2(\Omega)}^{4/5} \norm{(\cpm)^{p_k/2}}_{L^1(\Omega)}^{1/5} + \norm{(\cpm)^{p_k/2}}_{L^1(\Omega)}\Big)\norm{\nabla (\cpm)^{p_k/2}}_{L^2(\Omega)}\\
    &\quad = Cp_k \Big(\norm{\nabla (\cpm)^{p_k/2}}_{L^2(\Omega)}^{9/5} \norm{(\cpm)^{p_k/2}}_{L^1(\Omega)}^{1/5} + \norm{(\cpm)^{p_k/2}}_{L^1(\Omega)} \norm{\nabla (\cpm)^{p_k/2}}_{L^2(\Omega)}\Big)\\
    &\quad \leq \norm{\nabla (\cpm)^{p_k/2}}_{L^2(\Omega)}^2 + C ( p_k^{10} + p_k^2 )\norm{(\cpm)^{p_k/2}}_{L^1(\Omega)}^2  
    \\&\quad\leq
    \norm{\nabla (\cpm)^{p_k/2}}_{L^2(\Omega)}^2
    + C \bigl( p_k^{10} + 1 \bigr) \norm{(\cpm)^{p_k/2}}_{L^1(\Omega)}^2 .
\end{align*}
where we used Young's inequality twice, once with the exponents $(p, p') = (10/9, 10)$ and once with $(p, p') = (2,2)$.
Since $\frac{4 (p_k-1)}{p_k} = 4 - 4/p_k \geq 3$ for $k \geq 2$,
we can absorb the gradient term into the dissipation term on the left-hand side of~\eqref{eq:moser_iteration_1}.
Since $p_k/2 = p_{k-1}$ we arrive at
\begin{align}\label{eq:moser_iteration_2}
    \frac{\diff}{\diff t} \int_\Omega (\cpm)^{p_k} \diff \X
    \leq
        C \, (p_k^{10} + 1) \left( \int_\Omega (\cpm)^{p_{k-1}} \diff \X \right)^{2}
\end{align}
with a constant $C$ that depends on $\norm{\nabla \psi}_{L^\infty(0,T;L^6(\Omega))}$ but is independent of $k$.
The right-hand side of \eqref{eq:moser_iteration_2} no longer contains $\int_\Omega (\cpm)^{p_k} \diff \X$,
and its constant grows only polynomially in $p_k$.
Next, we integrate \eqref{eq:moser_iteration_2} in time and take the supremum over $t \in [0,T]$:
\begin{align*}
    A_k^{p_k}
    \leq |\Omega| \, \norm{\cpm_0}_{L^\infty(\Omega)}^{p_k} + C \, T \, (p_k^{10} + 1) \, A_{k-1}^{2 p_{k-1}}
    \leq C  (p_k^{10} + 1) \max\left(A_{k-1}, \norm{\cpm_0}_{L^\infty(\Omega)} \right)^{p_k} 
\end{align*}
where we also have absorbed $T$ and $|\Omega|$ in the constant.
Then, we take the $p_k$-th root of both sides of the inequality. Without loss of generality, we now assume that $C \geq 1$. Since $k \geq 2$ is arbitrary and on the right-hand side we have $A_{k-1}$, we can recursively apply this inequality. We obtain
\begin{align*}
    A_k
    &\leq
        (C(p_k^{10} + 1))^{1/p_k} \max\left(A_{k-1}, \norm{\cpm_0}_{L^\infty(\Omega)} \right)\\
    &\leq
        (C(p_k^{10} + 1))^{1/p_k} \max\left( (C(p_{k-1}^{10} + 1))^{1/p_{k-1}} \max\left(A_{k-2}, \norm{\cpm_0}_{L^\infty(\Omega)} \right), \norm{\cpm_0}_{L^\infty(\Omega)} \right)\\
    & =  (C(p_k^{10} + 1))^{1/p_k} \Big[ (C(p_{k-1}^{10} + 1))^{1/p_{k-1}} \max\left(A_{k-2}, \norm{\cpm_0}_{L^\infty(\Omega)} \right)\Big]\\
    & \leq \dots \leq 
    \prod_{j=2}^{k} (C (p_j^{10} + 1))^{1/p_j} \max\left(A_{1}, \norm{\cpm_0}_{L^\infty(\Omega)} \right),
    \end{align*}
where the equality in the middle follows from
\begin{equation*}
    (C(p_{k-1}^{10} + 1))^{1/p_{k-1}} \max\left(A_{k-2}, \norm{\cpm_0}_{L^\infty(\Omega)} \right) \geq  1 \cdot  \max\left(A_{k-2}, \norm{\cpm_0}_{L^\infty(\Omega)} \right) \geq \norm{\cpm_0}_{L^\infty(\Omega)}
\end{equation*}
 which holds because $C \geq 1$.\\
Next, we aim to pass to the limit $k \to \infty$.
To show that we get an $L^\infty$-bound for $\cpm$, we need to show that the product is convergent.
For that, we use the fact that an infinite product is convergent if and only if the sum of the logarithms of its factors is convergent. We have
\[
    \sum_{j=2}^{k} \frac{1}{p_j} \ln\bigl(C (p_j^{10} + 1) \bigr)
    = \sum_{j=2}^{k} \frac{\ln C + \ln (p_j^{10} + 1)}{2^{j}}
    \leq
    \sum_{j=2}^{k} \frac{\ln C + 10 j \ln (2)  + 1}{2^{j}}
    < \infty,
\]
where we used that $\ln(x + 1) \leq \ln(x) + 1$ for $x \geq 1$ and the convergence of the sum follows from the ratio test.
We can now estimate the $L^\infty$-bound for $\cpm$ by
\begin{align*}
    \norm{\cpm}_{L^\infty(Q)} &= \lim_{k \to \infty} \norm{\cpm}_{L^{p_k}(Q)} \leq \liminf_{k \to \infty} T^{\frac{1}{p_k}}\norm{\cpm}_{L^{\infty}(0,T;L^{p_k}(\Omega)}\\
    &= \liminf_{k \to \infty} T^{\frac{1}{p_k}} A_k 
    \leq C \max\left(A_1, \norm{\cpm_0}_{L^\infty(\Omega)} \right) \prod_{j=2}^{\infty} (C (p_j^{10} + 1))^{1/p_j}
    \leq C \max\left(A_1, \norm{\cpm_0}_{L^\infty(\Omega)}\right).
\end{align*}
This finishes the proof.

\subsection{Continuous dependence and uniqueness (proof of \zcref{thm:uniqueness})}
Let $(\bv_1, \cpm_1, \psi_1)$ and $(\bv_2, \cpm_2, \psi_2)$ be two weak solutions related, respectively, to the external potentials $\xi_1$, $\xi_2$ and to the initial data $(\bv_{0,1}, \cpm_{0,1})$, $(\bv_{0,2}, \cpm_{0,2})$. For brevity, we denote their differences as 
\begin{alignat*}{3}
    &\bv \coloneqq \bv_1 - \bv_2, \quad \quad && \cpm \coloneqq \cpm_1 - \cpm_2, \quad \quad && \psi \coloneqq \psi_1 - \psi_2\\
    &\bv_0 \coloneqq \bv_{0,1} - \bv_{0,2}, \quad \quad && \cpm_0 \coloneqq \cpm_{0,1} - \cpm_{0,2}, \quad \quad && \xi \coloneqq \xi_1 - \xi_2.
\end{alignat*}
Since during the proof we are going to need that $\cpm_1 \in L^{\infty}(Q)$ (and this regularity is not required by \zcref{definition:weak_sol}), let us assume that $(\bv_1, \cpm_1, \psi_1)$ is the solution constructed in \zcref{thm:weak_existence}, while $(\bv_2, \cpm_2, \psi_2)$ is an arbitrary weak solution. In the end, we will comment on it, explaining why this is not restrictive.
We test the difference of the equations with the difference of the corresponding unknowns, multiplying the $\psi$ equation by a big constant $K>0$ (fixed but yet to be chosen)---so that we can later on absorb $\psi$ terms from the Stokes and Nernst--Planck equations into the left-hand side---, and rewriting the nonlinearities by adding and subtracting terms. Summing everything, we have 
\begin{equation*}
    \begin{split}
        &\frac{1}{2} \frac{\diff}{\diff t}\bigg(\int_\Omega \Big(\kappa |\bv|^2 + \sigmp{|\cpm|^2}\Big) \diff \X \bigg) + \int_{\Omega} \Big( |\nabla \bv|^2 + \sigmp{|\nabla \cpm|_{\varMa}^2} + K |\nabla \psi|_{\varMa}^2\Big) \diff \X + \tau K \int_\Gamma |\psi|^2 \dS\\
        &\quad = - \int_\Omega (c^+_1 - c^-_1) \nabla \psi \cdot \bv \diff \X - \int_\Omega (c^+ - c^-) \nabla \psi_2 \cdot \bv \diff \X  + \int_\Omega \sigmp{\cpm_1  \bv \cdot \nabla \cpm } \diff \X \\
        &\quad \quad - \int_\Omega \sigmp{\pm \cpm_1 \varMa \nabla \psi \cdot \nabla \cpm \pm  \varMa \cpm \nabla \psi_2 \cdot \nabla \cpm} \diff \X + K \int_\Gamma \xi \psi \dS + K \int_\Omega (c^+ - c^-) \psi \diff \X.
    \end{split}
\end{equation*}
We bound the right-hand side by the boundedness of $\varMa$, the H\"older inequality, the fact that $\nabla\psi_2 \in L^\infty(0,T;L^6(\Omega))$ by \zcref{definition:weak_sol} (with a bound depending on $R$ only through \eqref{eq:main_norm_ineq} once uniqueness is established, see the end of the proof), and the boundedness of the charges \eqref{eq:bound_charges}. We obtain:
\begin{equation*}
    \begin{split}
        &\frac{1}{2} \frac{\diff}{\diff t}\bigg(\int_\Omega \Big(\kappa |\bv|^2 + \sigmp{|\cpm|^2}\Big) \diff \X \bigg) + \int_{\Omega} \Big( |\nabla \bv|^2 + \sigmp{|\nabla \cpm|_{\varMa}^2} + K |\nabla \psi|_{\varMa}^2\Big) \diff \X + \tau K \int_\Gamma |\psi|^2 \dS\\
        &\quad \leq \sigmp{\|\cpm_1\|_{L^{\infty}(Q)}} \|\nabla \psi\|_{L^2(\Omega)} \|\bv \|_{L^2(\Omega)} + \sigmp{\| \cpm\|_{L^3(\Omega)}} \| \nabla \psi_2 \|_{L^{\infty}(0,T;L^6(\Omega))} \| \bv \|_{L^2(\Omega)}\\
        &\quad \quad + \sigmp{ \| \cpm_1 \| _{L^{\infty}(Q)} \| \bv \|_{L^2(\Omega)} \| \nabla \cpm \|_{L^2(\Omega)}} + C \sigmp{\|\cpm_1\|_{L^{\infty}(Q)} \|\nabla \psi\|_{L^2(\Omega)} \|\nabla \cpm \|_{L^2(\Omega)} }\\
        &\quad\quad + C \! \sigmp{\| \cpm\|_{L^3(\Omega)} \| \nabla \psi_2 \|_{L^{\infty}(0,T;L^6(\Omega))} \| \nabla \cpm\|_{L^2(\Omega)}} + K\|\xi\|_{L^2(\Gamma)} \| \psi \|_{L^2(\Gamma)} + K\!\sigmp{\| \cpm \|_{L^2(\Omega)}} \|\psi\|_{L^2(\Omega)}\\
        &\quad \leq C \bigg(\|\nabla \psi\|_{L^2(\Omega)} \|\bv \|_{L^2(\Omega)} + \sigmp{\| \cpm\|_{L^3(\Omega)}} \| \bv \|_{L^2(\Omega)} + \| \bv \|_{L^2(\Omega)} \sigmp{\| \nabla \cpm \|_{L^2(\Omega)}} \\
        &\quad \quad + \|\nabla \psi\|_{L^2(\Omega)} \sigmp{ \|\nabla \cpm \|_{L^2(\Omega)} }  + \sigmp{\| \cpm\|_{L^3(\Omega)} \| \nabla \cpm\|_{L^2(\Omega)}}\bigg) + K \bigg( \|\xi\|_{L^2(\Gamma)} \| \psi \|_{L^2(\Gamma)}\\
        &\quad\quad + \sigmp{\| \cpm \|_{L^2(\Omega)}} \|\psi\|_{L^2(\Omega)}\bigg).
    \end{split}
\end{equation*}
Then, we use the Young inequality with a small constant $\delta > 0$ yet to be chosen, the Gagliardo--Nirenberg inequality \eqref{eq:GN_L3}, and the generalized Poincaré inequality \eqref{eq:gen_Poincare_ineq}, leading to:
    \begin{equation*}
    \begin{split}
        &\frac{1}{2} \frac{\diff}{\diff t}\bigg(\int_\Omega \Big(\kappa |\bv|^2 + \sigmp{|\cpm|^2}\Big) \diff \X \bigg) + \int_{\Omega} \Big( |\nabla \bv|^2 + \sigmp{|\nabla \cpm|_{\varMa}^2} + K |\nabla \psi|_{\varMa}^2\Big) \diff \X + \tau K \int_\Gamma |\psi|^2 \dS\\
        &\quad \leq \delta \Big( \sigmp{\| \nabla \cpm \|_{L^2(\Omega)}^2} + K  \| \psi \|_{L^2(\Gamma)}^2 + K \|\psi\|_{L^2(\Omega)}^2\Big)\\
        &\quad \quad + C_{\delta} \Big( \|\nabla \psi\|_{L^2(\Omega)}^2 + \|\bv \|_{L^2(\Omega)}^2 + \sigmp{\| \cpm\|_{L^3(\Omega)}^2} +  K \|\xi\|_{L^2(\Gamma)}^2  + K\sigmp{\| \cpm \|_{L^2(\Omega)}^2}\Big)\\
        &\quad \leq \delta \Big( \sigmp{\| \nabla \cpm \|_{L^2(\Omega)}^2} + K \| \psi \|_{L^2(\Gamma)}^2 + 2K C_P^2(\|\nabla \psi\|_{L^2(\Omega)}^2 + \|\psi\|_{L^2(\Gamma)}^2)\Big)\\
        &\quad \quad + C_{\delta} \Big( \|\nabla \psi\|_{L^2(\Omega)}^2 + \|\bv \|_{L^2(\Omega)}^2 +  K \|\xi\|_{L^2(\Gamma)}^2  + (K + 1)\sigmp{\| \cpm \|_{L^2(\Omega)}^2}\Big).
    \end{split}
\end{equation*}
Notice that every term multiplied by $\delta$ is either independent of $K$ or carries a factor $K$ that is matched by the corresponding term on the left-hand side ($K\|\nabla\psi\|^2_{\varMa}$ and $\tau K\|\psi\|^2_{L^2(\Gamma)}$). 
Hence we can choose $\delta>0$ sufficiently small, independently of $K$, to absorb all the terms multiplied by $\delta$ into the left-hand side. 
This fixes $C_\delta$, which depends only on the data and on $\delta$, but not on $K$. 
Subsequently, we choose $K > C_\delta$ large enough so that $C_\delta\|\nabla\psi\|^2_{L^2(\Omega)}$ is absorbed by the remaining $K(1-2\delta C_P^2) \|\nabla\psi\|^2_{L^2(\Omega)}$. 
Finally, we integrate in time on $(0,s)$ for $s \in (0,T)$, obtaining:
\begin{equation*}
    \begin{split}
       &\int_\Omega \Big(\kappa |\bv|^2 + \sigmp{|\cpm|^2}\Big) \diff \X + \int_0^s\int_{\Omega} \Big( |\nabla \bv|^2 + \sigmp{|\nabla \cpm|_{\varMa}^2} + |\nabla \psi|_{\varMa}^2\Big) \diff \X \diff t + \tau \int_0^s\int_\Gamma |\psi|^2 \dS \diff t\\
       & \quad \leq  \kappa \| \bv_0 \|_{L^2(\Omega)}^2 +  \sigmp{\| \cpm_0 \|_{L^2(\Omega)}^2} + C \int_0^s \|\xi\|_{L^2(\Gamma)}^2 \diff t + C \int_0^s \Big( \| \bv \|_{L^2(\Omega)}^2 + \sigmp{\| \cpm \|_{L^2(\Omega)}^2} \Big) \diff t.
    \end{split}
\end{equation*}
By the Gronwall inequality, we conclude that:
\begin{equation*}
    \begin{split}
        &\| \bv \|_{L^{\infty}(0,T;L^2(\Omega)) \cap L^2(0,T;W^{1,2}(\Omega))}^2 + \sigmp{\| \cpm\|_{L^{\infty}(0,T;L^2(\Omega)) \cap L^2(0,T;W^{1,2}(\Omega))}^2 } +  \| \psi \|_{ L^2(0,T;W^{1,2}(\Omega))}^2\\
        &\quad \leq C \Big( \| \bv_0 \|_{L^2(\Omega)}^2 +  \sigmp{\| \cpm_0 \|_{L^2(\Omega)}^2} + \| \xi\|_{L^2(0,T;L^2(\Gamma))}^2\Big),
    \end{split}
\end{equation*}
where we have also used the generalized Poincaré inequality \eqref{eq:gen_Poincare_ineq} for $\psi$.
Improving the estimate for the electric potential can be easily done as we did at the end of Step \ref{enum:grad_est_charges}, since the equation \eqref{eq:weak_Poisson} is linear. Invoking \eqref{eq:est_psii} for $p=2$ (together with the embedding $W^{2,2}(\Omega) \hookrightarrow W^{1,6}(\Omega)$) and \eqref{eq:est_psie} for $q=\infty$, $p=6$,  we get
\begin{equation*}
    \begin{split}
        \norm{\psi}_{L^{\infty}(0,T;W^{1,6}(\Omega))}^2 &\leq C \Big( \norm{c^+ - c^-}_{L^{\infty}(0,T;L^2(\Omega))}^2 + \norm{\xi}_{L^{\infty}(\Sigma)}^2\Big)\\
        &\leq C \Big( \| \bv_0 \|_{L^2(\Omega)}^2 +  \sigmp{\| \cpm_0 \|_{L^2(\Omega)}^2} + \norm{\xi}_{L^{\infty}(\Sigma)}^2 \Big),
    \end{split}
\end{equation*}
where the second inequality follows from the estimate above, used to bound the norm of the charges.
Finally, uniqueness is an easy consequence of \eqref{eq:cont_dep_ineq}: taking two solutions with the same initial data and the same external electric potential, the first one constructed as in \zcref{thm:weak_existence} and the second one arbitrary, the right-hand side vanishes. Thus, the two solutions must coincide a.e., proving uniqueness of the solution. In particular, every weak solution satisfies the estimates \eqref{eq:main_norm_ineq} and \eqref{eq:bound_charges}, so that \eqref{eq:cont_dep_ineq} holds for every couple of weak solutions.

\section{The optimal control problem}\label{sec:optimal_control}
Thanks to \zcref{thm:weak_existence} and \zcref{thm:uniqueness}, for every control $\xi \in \CXR$ there exists a unique associated state $(\bv,\cpm,\psi)$ solving system \eqref{eq:snpp}--\eqref{eq:ic} in the sense of \zcref{definition:weak_sol}. Thus, the following \textit {solution operator}---also called \textit{control-to-state operator}---is well-defined: 
\begin{equation*}
    \Sol \colon \CXR \subset L^{\infty}(\Sigma) \to \CY, \qquad \xi \mapsto (\bv,\cpm,\psi),
\end{equation*}
where the Banach space $\CY$ is 
\begin{equation*}
    \begin{split}
        \CY \coloneqq  &\left[W^{1,2}(0,T;L^2(\Omega)) \cap L^2(0,T;W^{2,2}(\Omega)) \cap L^2(0,T;W^{1,2}_{0,\sigma}(\Omega))\right]\\
        &\times \left[W^{1,2}(0,T;(W^{1,2}(\Omega))^*) \cap L^\infty(0,T;L^2(\Omega)) \cap L^2(0,T;W^{1,2}(\Omega)) \cap L^{\infty}(Q)\right]\\
    &\times L^\infty(0,T;W^{1,6}(\Omega)).
    \end{split}
\end{equation*}
Moreover, introducing the larger Banach space
\begin{equation*}
    \begin{split}
        \CZ \coloneqq  &\left[L^{\infty}(0,T;L^2(\Omega)) \cap L^2(0,T;W^{1,2}(\Omega))\right]\\
        &\times \left[L^{\infty}(0,T;L^2(\Omega)) \cap L^2(0,T;W^{1,2}(\Omega)) \right]\\
        &\times L^{\infty}(0,T;W^{1,6}(\Omega)),
    \end{split}
\end{equation*}
the continuous dependence inequalities \eqref{eq:cont_dep_ineq}--\eqref{eq:improved_cont_dep_ineq} proved in \zcref{thm:uniqueness} guarantee that $\Sol \colon \CXR \subseteq L^{\infty}(\Sigma) \to \CZ$ is Lipschitz continuous.
We can then introduce the \textit{reduced cost functional}
\begin{equation*}
    \redJ(\xi) \coloneqq \J(\xi, \Sol(\xi)), \qquad \xi \in \CXR,
\end{equation*}
and, with this notation, state the optimal control problem under examination more compactly as
\begin{equation}\label{eq:opt_control_prob}
    \min_{\xi \in \CXad}{\redJ}(\xi). 
\end{equation}
We first aim to prove that the minimization problem \eqref{eq:opt_control_prob} admits at least one solution, henceforth called an \textit{optimal control}.

\begin{thm}[Existence of an optimal control]
   There exists at least one optimal control $\xi^* \in \CXad$ solving problem \eqref{eq:opt_control_prob}.
\end{thm}
\begin{proof}
    This follows directly from the estimate \eqref{eq:main_norm_ineq} in \zcref{thm:weak_existence} and the direct method of the Calculus of Variations.
    First of all, we notice that the cost functional is bounded from below since it is the sum of
    \begin{equation*}
        -\frac{\alpv }{2}\int_0^T\int_\Omega |\nabla \times \bv|^2 \diff \X \diff t \geq - C \|\bv\|_{L^2(0,T;W^{1,2}(\Omega))}^2 \geq - C C_{1,R}^2, 
    \end{equation*}
    and three nonnegative terms. Thus, $\inf_{\xi \in \CXad} \redJ(\xi)$ is finite.
    Let $(\xi_n)_n$ be a minimizing sequence, i.e., 
    \[
        \lim_{n \to \infty} \redJ(\xi_n) = \inf_{\xi \in \CXad} \redJ(\xi),
    \]
    and denote $(\bv_n,\cpm_n,\psi_n) \coloneqq \Sol(\xi_n)$ the corresponding sequence of states. 
    Since $\CXad$ is bounded and weakly-$\ast$ closed in $L^{\infty}(\Sigma)$, there exists a $\xi^* \in \CXad$ such that, at least along a subsequence that we do not relabel, 
    \begin{alignat}{2}\label{eq:xi_weak_conv}
        \xi_n \overset{\ast}{\rightharpoonup} \xi^{*} &\quad \text{in } L^{\infty}(\Sigma).
    \end{alignat}
    Moreover, by the estimate \eqref{eq:main_norm_ineq}, there exists $(\bv^*,{\cpm}^*,\psi^*) \in \CY$ such that, along a nonrelabelled subsequence, the following convergences hold by Banach--Alaoglu (see, e.g., \cite[Thm.~3.16]{brezis_11}) and Aubin--Lions (see, e.g., \cite[Lem.~7.7]{roubicek_13}) Theorems
    \begin{alignat}{2}
        \bv_n &\rightharpoonup \bv^{*} &&\quad \text{in } W^{1,2}(0,T;L^2(\Omega)) \cap L^2(0,T;W^{1,2}_{0,\sigma}(\Omega) \cap W^{2,2}(\Omega)), \label{eq:v_weak_conv}\\
        \bv_n &\to \bv^{*} &&\quad \text{in } L^2(0,T;W^{1,2}_{0,\sigma}(\Omega)),\\
        \cpm_n &\overset{\ast}{\rightharpoonup} {\cpm}^{*} &&\quad \text{in } W^{1,2}(0,T;(W^{1,2}(\Omega))^*) \cap L^\infty(0,T;L^2(\Omega)) \cap L^2(0,T;W^{1,2}(\Omega)),\label{eq:c_weak_conv}\\
        \cpm_n &\to {\cpm}^{*} &&\quad \text{in } L^2(0,T;L^2(\Omega)),\\
        \cpm_n &\to {\cpm}^{*} &&\quad \text{a.e. } \text{in } Q, \label{eq:c_ae_conv}\\
        \psi_n &\overset{\ast}{\rightharpoonup} \psi^{*} &&\quad \text{in }L^\infty(0,T;W^{1,6}(\Omega)). \label{eq:psi_weak_conv}
    \end{alignat}
    Moreover, the boundedness of the charges in $L^\infty(Q)$, the pointwise convergence \eqref{eq:c_ae_conv}, and the Dominated Convergence Theorem yield 
    \begin{alignat}{2}
        \cpm_n \to {\cpm}^{*} &\quad \text{in } L^p(0,T;L^p(\Omega)) \label{eq:c_sconv_p}
    \end{alignat}
    for all $p < +\infty$. Then, the convergences \eqref{eq:v_weak_conv}, \eqref{eq:psi_weak_conv}, and \eqref{eq:c_sconv_p} ensure that 
     \begin{alignat}{2}
        \cpm_n \nabla \psi_n \rightharpoonup {\cpm}^* \nabla \psi^* \quad \text{and} \quad \cpm_n \bv_n \rightharpoonup {\cpm}^* \bv^* &\quad \text{in } L^2(0,T;L^2(\Omega)).\label{eq:product_weak_conv}
    \end{alignat}
    The convergences \eqref{eq:xi_weak_conv}--\eqref{eq:product_weak_conv} suffice to pass to the limit in the weak formulation \eqref{eq:weak_snpp} written for $\xi_n$ and we get that $(\bv^*, {\cpm}^*, \psi^*)$ is a solution of the same system written for $\xi^*$. 
    Finally, $\bv^*$ and ${\cpm}^*$ satisfy the initial conditions \eqref{eq:ic} in $L^2(\Omega)$. 
    Indeed, by \eqref{eq:v_weak_conv} it follows that $\bv_n(0)=\bv_0 \rightharpoonup \bv^*(0)$ in $L^2(\Omega)$, thus $\bv^*(0) = \bv_0$ in $L^2(\Omega)$. By \eqref{eq:c_weak_conv}, we deduce the weaker $\cpm_n(0)=\cpm_0 \rightharpoonup {\cpm}^*(0)$ in $(W^{1,2}(\Omega))^*$, from which ${\cpm}^*(0) = \cpm_0$ in $(W^{1,2}(\Omega))^*$ (that can be extended to $L^2(\Omega)$ since ${\cpm}^* \in C^0([0,T];L^2(\Omega))$).
    Thus, $(\bv^{*}, {\cpm}^{*}, \psi^{*})$ is a weak solution to the SNPP system \eqref{eq:snpp}--\eqref{eq:ic} with control $\xi^*$, i.e., that $(\bv^*, {\cpm}^*, \psi^*) = \Sol(\xi^*)$. 
    As the cost functional is not weakly lower semicontinuous due to the minus sign in front of the curl of the velocity field, we need the strong convergence of the velocity field to be able to pass to the limit in the cost functional:
    \begin{align*}
        &\inf_{\xi \in \CXad} \redJ(\xi)= \lim_{n \to \infty} \redJ(\xi_n)\\
        &\quad  = \lim_{n \to +\infty} \left(\int_0^T \!\int_\Omega \Big( - \frac{\alpv }{2}|\nabla \times \bv_n|^2 + \frac{1}{2} \sigmp{\alpcpm |c_n^{\pm} - c^{\pm}_Q|^2}\Big) \diff \X \diff t + \frac{\alpxi}{2} \int_0^T \! \int_\Gamma |\xi_n|^2 \dS \diff t\right)\\
        &\quad \geq
           \lim_{n \to +\infty} \! \left(\int_0^T \!\int_\Omega \Big(- \frac{\alpv }{2} |\nabla \times \bv_n|^2 + \frac{1}{2} \sigmp{\alpcpm |c_n^{\pm} - c^{\pm}_Q|^2}  \Big)\diff \X \diff t \right) +  \liminf_{n \to +\infty} \left(\frac{\alpxi}{2} \int_0^T \! \int_\Gamma |\xi_n|^2 \dS \diff t \right)\\
        &\quad \geq
            \int_0^T \int_\Omega  \Big(- \frac{\alpv }{2} |\nabla \times \bv^*|^2 + \frac{1}{2} \sigmp{\alpcpm |{c^{\pm}}^{*} - c^{\pm}_Q|^2}  \Big)\diff \X \diff t + \frac{\alpxi}{2} \int_0^T \int_\Gamma |\xi^*|^2 \dS \diff t 
        = \redJ(\xi^*),
    \end{align*}
    where we have used the strong convergence of the velocity field and the charges for the convergence of the first term and the weak-$\ast$ convergence of $\xi_n$ and the weak lower semicontinuity of the norm for the lower bound of the second term.
    From this inequality, the optimality of the limit $\xi^*$ follows.
\end{proof}

\noindent The subsequent more challenging goal consists in characterizing the optimal controls through some \textit{first-order necessary optimality conditions}, which will have the form of a variational inequality (see \zcref{thm:opt_cond} below). 
To derive them, we consider an optimal control $\xi^* \in \CXad$ and any other admissible control $\xi \in \CXad$. Since $\CXad$ is convex, for every $\theta \in (0,1)$ their convex combination $\xi^* + \theta(\xi-\xi^*)$ is still an admissible control. Thus, by the minimality of $\xi^*$, we have:
\begin{equation*}
    \frac{\redJ\big(\xi^* + \theta(\xi-\xi^*)\big) - \redJ(\xi^*)}{\theta}\geq 0
\end{equation*}
for every $\theta \in (0,1)$ and for every $\xi \in \CXad$. Knowing that $\redJ$ is at least G\^{a}teaux differentiable, it would be possible to pass to the limit as $\theta \to 0^+$ in this inequality, showing that
\begin{equation}\label{eq:opt_cond_0}
    D\redJ(\xi^*)[\xi -\xi^*] \geq 0 \qquad \forall \xi \in \CXad,
\end{equation}
which are the conditions we are looking for. Thus, to rigorously prove that \eqref{eq:opt_cond_0} holds, we just need to prove that $\redJ(\cdot) = \J(\cdot,\Sol(\cdot))$ is at least G\^{a}teaux differentiable.  Since the reduced cost functional is the composition of the cost functional $\J$ and the solution operator $\Sol$, first of all we notice that $\J$ is trivially Fréchet differentiable on $\CXR \times \CZ$. In fact, it is the algebraic combination of the squares of $L^2$-norms, and $\CXR \times \CZ$ is an open subset of
\begin{equation*}
    L^{\infty}(\Sigma) \times \CZ  \, \, \hookrightarrow \, \, L^2(\Sigma) \times \Big[ L^2(0,T;W^{1,2}(\Omega)) \times \big(L^2(Q)\big)^2 \Big].
\end{equation*}
In the next section, we are going to prove that also $\Sol$ is Fréchet differentiable as an operator between suitable Banach spaces. Thus, by the chain rule, their composition $\redJ$ is Fréchet differentiable and consequently G\^{a}teaux differentiable. 

\section{Differentiability of the control-to-state operator}\label{sec:differentiability}

The Fr\'{e}chet derivative of the control-to-state operator $\Sol$ at the point $\xi$ in direction $\lxi$ is classically given by the solution of the linearized SNPP system $(\lv, \lcpm, \lpsi)$ at the point $(\xi,\Sol(\xi))$, which, written in the strong form, is:
\begin{subequations}\label{eq:lin_snpp}
    \begin{align}
        \kappa \, \partial_t \lv - \Delta \lv + \nabla \lp + (\lc^+ - \lc^-) \nabla \psi + (c^+ - c^-) \nabla \lpsi &= 0 &&\text{ in } Q,\label{eq:lin_s}\\
        \nabla \cdot \lv &= 0 &&\text{ in } Q,\label{eq:lin_incomp}\\
        \partial_t \lcpm + \nabla \cdot (\lcpm \bv + \cpm \lv) - \nabla \cdot (\varMa (\nabla \lcpm \pm \lcpm \nabla \psi \pm \cpm \nabla \lpsi)) &= 0 &&\text{ in } Q,\label{eq:lin_np}\\
        -\nabla \cdot (\varMa \nabla \lpsi) - (\lc^+ - \lc^-) &= 0 &&\text{ in } Q,\label{eq:lin_p}
\intertext{equipped with the boundary conditions\newsubeqgroup{eq:lin_bc}}
        \lv &= 0 &&\text{ on } \Sigma,\\
        \varMa (\nabla \lcpm \pm \cpm \nabla \lpsi \pm \lcpm \nabla \psi) \cdot \bn &= 0 &&\text{ on } \Sigma,\\
        \varMa \nabla \lpsi \cdot \bn + \tau \lpsi &= \lxi  &&\text{ on } \Sigma,
\intertext{and the initial conditions\stopsubeqgroup}
        \lv(0) = 0, \quad \lcpm(0) &= 0 &&\text{ in } \Omega. \label{eq:lin_ic}
    \end{align}
\end{subequations}

\begin{remark}
    Notice that, for a fixed control $\xi \in \CXR$, we have an infinite number of linearized systems, one for each $\lxi$. However, to keep the notation as simple as possible, we denote the solution to \eqref{eq:lin_snpp}--\eqref{eq:lin_ic}  corresponding to a source term $\lxi$ simply by $(\lv, \lcpm,\lpsi)$, suppressing the dependence on $\lxi$.
\end{remark}

\noindent First of all, we introduce our definition of a weak solution.

\begin{definition}[Weak solution of the linearized system]\label{definition:lin_weak_sol}
    Let $(\bv, \cpm, \psi) \in \CY$ and let
    $\lxi \in L^2(0,T;L^2(\Gamma))$.  A tuple $(\lv,\lcpm,\lpsi)$ with
    \begin{align*}
      \lv &\in W^{1,2}(0,T;L^2(\Omega))\cap L^2(0,T;W^{1,2}_{0,\sigma}(\Omega) \cap W^{2,2}(\Omega)),\\
      \lcpm &\in W^{1,2}(0,T;(W^{1,2}(\Omega))^*)\cap C^0([0,T];L^2(\Omega))\cap L^2(0,T;W^{1,2}(\Omega)),\\
      \lpsi &\in L^2(0,T;W^{1,2}(\Omega))
    \end{align*}
    is called a weak solution of \eqref{eq:lin_snpp}--\eqref{eq:lin_ic} if the system \eqref{eq:lin_snpp}--\eqref{eq:lin_bc} is satisfied in the following variational sense
    \begin{subequations} \label{eq:lin_weak}
        \begin{align}
              \kappa\,\int_\Omega \partial_t\lv \cdot \bs{\varphi} \diff \X +\int_\Omega\nabla\lv:\nabla\bs{\varphi}\diff \X +\int_\Omega\Big[(\gamma^+ - \gamma^-) \nabla\psi + ( c^+- c^-)\nabla\lpsi\Big]\cdot\bs{\varphi}\diff \X&=0, \label{eq:lin_weak_s}\\
              \langle\partial_t\lcpm,\varphi\rangle - \int_\Omega\big(\lcpm \bv+ c^{\pm}\lv\big)\cdot\nabla\varphi\diff \X + \int_\Omega\varMa \big(\nabla\lcpm\pm\lcpm\nabla\psi \pm c^{\pm}\nabla\lpsi\big)\cdot\nabla\varphi\diff \X&=0, \label{eq:lin_weak_np}\\
              \int_\Omega \varMa \nabla\lpsi\cdot\nabla\varphi\diff \X +\int_\Gamma\big(\tau\lpsi-\lxi\big)\varphi \dS -\int_\Omega(\gamma^+-\gamma^-)\varphi\diff \X&=0, \label{eq:lin_weak_p}
        \end{align}
    \end{subequations}
    for almost every $t\in(0,T)$ and for all $\bs{\varphi}\in W^{1,2}_{0,\sigma}(\Omega)$ and all $\varphi\in W^{1,2}(\Omega)$, 
    and if the initial conditions \eqref{eq:lin_ic} hold in $L^2(\Omega)$.
\end{definition}

\begin{remark}\label{remark:pi}
    The pressure $\pi$ does not appear in the definition of weak solution because it can always be reconstructed a posteriori (cf. \zcref{remark:p} for the state system). Also in this case, by comparison in \eqref{eq:lin_s} and choosing $\pi$ having vanishing mean value, $\pi \in L^2(0,T;W^{1,2}(\Omega))$.
\end{remark}

\noindent Then, we need to prove that this system is well-posed.

\begin{propo}[Well-posedness of the linearized system]\label{propo:exist_lin_snpp}
    For every control $\xi \in \CXR$ with associated state $(\bv, \cpm, \psi) \coloneqq \Sol (\xi)$, and for every $\lxi \in L^2(0,T;L^2(\Gamma))$, the linearized system \eqref{eq:lin_snpp} has a unique weak solution $(\lv, \lcpm, \lpsi)$ in the sense of \zcref{definition:lin_weak_sol}.  Moreover, there exists a constant $C_{4,R}>0$ (independent of $\lxi$) such that the following inequality holds:
    \begin{equation}\label{eq:lin_syst_main_est}
        \begin{split}
            &\norm{\lv}_{W^{1,2}(0,T;L^2(\Omega))\cap L^\infty(0,T;L^2(\Omega)) \cap L^2(0,T; W^{2,2}(\Omega))}\\
            &\quad \quad + \norm{\lcpm}_{W^{1,2}(0,T;(W^{1,2}(\Omega))^*)\cap L^\infty(0,T;L^2(\Omega)) \cap L^2(0,T;W^{1,2}(\Omega))} + \norm{\lpsi}_{L^2(0,T;W^{1,2}(\Omega))}\\
            &\quad \leq C_{4,R} \norm{\lxi}_{L^2(0,T;L^2(\Gamma))}.
        \end{split}
    \end{equation}
    Finally, there exists another constant $C_{5,R} > 0$, such that, for every $\lxi \in L^{\infty}(\Sigma)$, it holds
    \begin{equation}\label{eq:lin_lpsi_add_est}
        \norm{\lpsi}_{L^{\infty}(0,T;W^{1,6}(\Omega))} \leq C_{5,R}\norm{\lxi}_{L^{\infty}(\Sigma)}.
    \end{equation}
\end{propo}

\begin{remark}
    In particular, we will prove that the solution of the linearized system $(\lv, \lcpm, \lpsi)$ for $h \in L^\infty(\Sigma)$ belongs to the Banach space $\CZ$.
\end{remark}

\begin{proof}
    First of all, we notice that uniqueness follows from the linearity of the system and from the estimate \eqref{eq:lin_syst_main_est}. 
    To prove existence, one should introduce a suitable discretization, derive some uniform a priori estimate, and pass to the limit. Here we will proceed formally, deriving the a priori estimates we need at the continuous level.
    
    \paragraph{Main a priori estimate.} Summing \eqref{eq:lin_weak_s} tested with $\bs{\varphi} = \lv$, \eqref{eq:lin_weak_np} tested with $\varphi = \lcpm$, and \eqref{eq:lin_weak_p} tested with $\varphi = \lpsi$ and again multiplied by a constant $K>0$ yet to be chosen, we have: 
        \begin{align*}
            &\frac{1}{2} \frac{\diff}{\diff t} \int_{\Omega}\left( \kappa |\lv|^2 + \sigmp{|\lcpm|^2}\right) \diff \X
            +  \int_\Omega \Big(|\nabla \lv|^2 + \sigmp{|\nabla \lcpm|^2_{\varMa}} + K |\nabla \lpsi|_{\varMa}^2\Big) \diff \X + \tau K \int_{\Gamma} |\lpsi|^2 \dS\\
            &\quad = - \int_{\Omega} \Big[(\gamma^+ - \gamma^-) \nabla\psi + ( c^+- c^-)\nabla\lpsi\Big]\cdot\lv \diff \X + \int_{\Omega} \sigmp{\cpm \lv \cdot \nabla \lcpm} \diff \X\\
            &\quad \quad - \int_\Omega\varMa \sigmp{\big(\pm\lcpm\nabla\psi \pm c^{\pm}\nabla\lpsi \big)\cdot\nabla\lcpm} \diff \X + K \int_{\Gamma} \lxi \lpsi \dS + K \int_{\Omega} (\gamma^+ - \gamma^-)\lpsi \diff \X.
        \end{align*}
    Then, we proceed to bound the right-hand side, using the H\"older inequalities and the estimates \eqref{eq:main_norm_ineq}--\eqref{eq:bound_charges}:
        \begin{align*}
            &\frac{1}{2} \frac{\diff}{\diff t} \int_{\Omega}\left( \kappa |\lv|^2 + \sigmp{|\lcpm|^2}\right) \diff \X
            +  \int_\Omega \Big(|\nabla \lv|^2 + \sigmp{|\nabla \lcpm|^2_{\varMa}} + K |\nabla \lpsi|_{\varMa}^2\Big) \diff \X + \tau K \int_{\Gamma} |\lpsi|^2 \dS\\
            &\quad \leq \sigmp{\norm{\lcpm}_{L^2(\Omega)}} \norm{\nabla\psi}_{L^{\infty}(0,T;L^6(\Omega))} \norm{\lv}_{L^3(\Omega)}+ \sigmp{\norm{\cpm}_{L^{\infty}(Q)}} \norm{\nabla\lpsi}_{L^2(\Omega)}\norm{\lv}_{L^2(\Omega)}\\
            &\quad \quad + \sigmp{\norm{\cpm}_{L^{\infty}(Q)}\norm{\lv}_{L^2(\Omega)} \norm{\nabla \lcpm}_{L^2(\Omega)}} 
            + \sigmp{\norm{\lcpm}_{L^3(\Omega)} \norm{\nabla \psi}_{L^{\infty}(0,T;L^6(\Omega))} \norm{\nabla \lcpm}_{L^2(\Omega)}}\\
            &\quad\quad+ \sigmp{\norm{\cpm}_{L^{\infty}(Q)} \norm{\nabla\lpsi}_{L^2(\Omega)} \norm{\nabla \lcpm}_{L^2(\Omega)}}
            + K \norm{\lxi}_{L^2(\Gamma)}\norm{\lpsi}_{L^2(\Gamma)} + K \sigmp{\norm{\lcpm}_{L^2(\Omega)}}\norm{\lpsi}_{L^2(\Omega)}\\
            &\quad \leq C_R \Big( \sigmp{\norm{\lcpm}_{L^2(\Omega)}} \norm{\lv}_{L^3(\Omega)}
            + \norm{\nabla\lpsi}_{L^2(\Omega)}\norm{\lv}_{L^2(\Omega)} 
            + \norm{\lv}_{L^2(\Omega)} \sigmp{\norm{\nabla \lcpm}_{L^2(\Omega)}}\\
            &\quad \quad + \sigmp{\norm{\lcpm}_{L^3(\Omega)}\norm{\nabla \lcpm}_{L^2(\Omega)}}
            + \norm{\nabla\lpsi}_{L^2(\Omega)} \sigmp{\norm{\nabla \lcpm}_{L^2(\Omega)}}\Big)\\
            &\quad \quad + K \norm{\lxi}_{L^2(\Gamma)}\norm{\lpsi}_{L^2(\Gamma)}
            + K \sigmp{\norm{\lcpm}_{L^2(\Omega)}}\norm{\lpsi}_{L^2(\Omega)}.
    \end{align*}
    Next, we apply the Gagliardo--Nirenberg inequality \eqref{eq:GN_L3} and the Young inequality with a small $\delta>0$, obtaining:
        \begin{align*}
        &\frac{1}{2} \frac{\diff}{\diff t} \int_{\Omega}\left( \kappa |\lv|^2 + \sigmp{|\lcpm|^2}\right) \diff \X
        +  \int_\Omega \Big(|\nabla \lv|^2 + \sigmp{|\nabla \lcpm|^2_{\varMa}} + K |\nabla \lpsi|_{\varMa}^2\Big) \diff \X + \tau K \int_{\Gamma} |\lpsi|^2 \dS\\
        &\quad \leq \delta \Big( \norm{\nabla \lv}_{L^2(\Omega)}^2 + \sigmp{\norm{\nabla \lcpm}_{L^2(\Omega)}^2} + K \norm{\lpsi}_{L^2(\Gamma)}^2  + K \norm{\lpsi}_{L^2(\Omega)}^2 \Big)\\
        &\quad \quad + C_{\delta}\Big( (1 + K)\sigmp{\norm{\lcpm}_{L^2(\Omega)}^2} + \norm{\lv}_{L^2(\Omega)}^2 + \norm{\nabla \lpsi}_{L^2(\Omega)}^2 + K \norm{\lxi}^{2}_{L^2(\Gamma)}\Big).
    \end{align*}
    Then,  we use the generalized Poincaré inequality \eqref{eq:gen_Poincare_ineq} to bound the term $\norm{\lpsi}_{L^2(\Omega)}$ on the right-hand side. Fixing $\delta$ small enough, the terms multiplied by it can be absorbed into the left-hand side. After that, choosing $K$ big enough, also $\norm{\nabla \lpsi}_{L^2(\Omega)}^2$ can be absorbed by the corresponding term on the left-hand side. Integrating in time and applying the Gronwall inequality, we conclude:
    \begin{equation}\label{eq:lin_syst_main_est_1}
        \begin{split}
                &\norm{\lv}_{L^\infty(0,T;L^2(\Omega)) \cap L^2(0,T;W^{1,2}(\Omega))}
                + \norm{\lcpm}_{L^\infty(0,T;L^2(\Omega)) \cap L^2(0,T;W^{1,2}(\Omega))}
                + \norm{\lpsi}_{L^2(0,T;W^{1,2}(\Omega))}\\
                &\quad \leq C_R \norm{\lxi}_{L^2(0,T;L^2(\Gamma))}.
        \end{split}
    \end{equation}
    
    \paragraph{Maximal regularity for $\lv$.} Rewriting \eqref{eq:lin_s} as
    \begin{equation*}
        \kappa \, \partial_t \lv - \Delta \lv + \nabla \lp = - (\lc^+ - \lc^-) \nabla \psi - (c^+ - c^-) \nabla \lpsi \eqqcolon f, 
    \end{equation*}
    we notice that 
    \begin{equation*}
        \begin{split}
            \norm{f}_{L^2(0,T;L^2(\Omega))} &\leq \sigmp{\norm{\lcpm}_{L^2(0,T;L^3(\Omega))}}\norm{\nabla \psi}_{L^{\infty}(0,T;L^6(\Omega))} + \sigmp{\norm{\cpm}_{L^{\infty}(Q)}}\norm{\nabla \lpsi}_{L^2(0,T;L^2(\Omega))}\\
            &\leq C_R  \norm{\lxi}_{L^2(0,T;L^2(\Gamma))}
        \end{split}
    \end{equation*}
    where we have employed the estimates \eqref{eq:main_norm_ineq} and \eqref{eq:lin_syst_main_est_1}. Thus, by maximal parabolic regularity for the nonstationary Stokes equation (see, e.g., \cite[Section IV, Thm.~2.7.1]{sohr_01}), we have:
    \begin{equation}\label{eq:lin_syst_main_est_2}
        \begin{split}
            \norm{\lv}_{W^{1,2}(0,T;L^2(\Omega))\cap L^2(0,T; W^{2,2}(\Omega))} \leq C \norm{f}_{L^2(0,T;L^2(\Omega))} \leq C_R \norm{\lxi}_{L^2(0,T;L^2(\Gamma))}.
        \end{split}
    \end{equation}
    
    \paragraph{Time regularity for $\lcpm$.} The a priori estimate for $\partial_t \lcpm$ can be derived by comparison in \eqref{eq:lin_np}, since
    \begin{equation*}
        \begin{split}
            &\norm{\lcpm \bv + \cpm \lv - \varMa ( \nabla \lcpm \pm \lcpm \nabla \psi \pm \cpm \nabla \lpsi)}_{L^2(0,T;L^2(\Omega))}\\
            &\quad\leq \norm{\lcpm}_{L^{\infty}(0,T;L^2(\Omega))} \norm{\bv}_{L^2(0,T;L^{\infty}(\Omega))} + \norm{\cpm}_{L^{\infty}(Q)}\norm{\lv}_{L^2(0,T;L^2(\Omega))}\\
            &\quad \quad + C \Big(\!\norm{\nabla \lcpm}_{L^2(0,T;L^2(\Omega))} \! + \norm{\lcpm}_{L^2(0,T;L^3(\Omega))} \! \norm{\nabla \psi}_{L^{\infty}(0,T;L^6(\Omega))} \!+ \norm{\cpm}_{L^{\infty}(Q)} \! \norm{\nabla \lpsi}_{L^2(0,T;L^2(\Omega))} \!\Big)\\
            &\quad \leq C_R \norm{\lxi}_{L^2(0,T;L^2(\Gamma))},
        \end{split}
    \end{equation*}
    where the last inequality follows again from \eqref{eq:main_norm_ineq} and \eqref{eq:lin_syst_main_est_1}. Thus, $\norm{\partial_t \lcpm}_{L^2(0,T; (W^{1,2}(\Omega))^*)} \leq C_R \norm{\lxi}_{L^2(0,T;L^2(\Gamma))}$.
    
    \paragraph{More regularity for $\lpsi$.} Since $\lc^+ - \lc^-$ is only bounded in $L^{\infty}(0,T;L^2(\Omega))$, we split $\lpsi = \lpsi_i + \lpsi_e$ as in \eqref{eq:splitting}, with $\lc^+ - \lc^-$ and $\lxi$ in place of $c^+ - c^-$ and $\xi$. Then \eqref{eq:est_psii} with $p = 2$, the embedding $W^{2,2}(\Omega) \hookrightarrow W^{1,6}(\Omega)$, \eqref{eq:lin_syst_main_est_1}, and \zcref{propo:gradient_integrability} with $F = 0$, $g = \lxi$, $p = 6$, $q = \infty$ yield
    \begin{equation*}
        \begin{split}
            \norm{\lpsi}_{L^{\infty}(0,T;W^{1,6}(\Omega))}
            &\leq C \Big( \norm{\lc^+ - \lc^-}_{L^{\infty}(0,T;L^2(\Omega))} + \norm{\lxi}_{L^{\infty}(\Sigma)} \Big)\\
            &\leq C_R \Big( \norm{\lxi}_{L^2(0,T;L^2(\Gamma))} + \norm{\lxi}_{L^{\infty}(\Sigma)} \Big)
            \leq C_R \norm{\lxi}_{L^{\infty}(\Sigma)},
        \end{split}
    \end{equation*}
    and the proof is complete.
\end{proof}

\noindent We are now in the position to prove what we have already claimed, which is that the solution of the linearized system is the Fréchet derivative of the solution operator.

\begin{propo}[Differentiability of the solution operator]
    The solution operator $\Sol: \CXR \subseteq L^{\infty}(\Sigma) \to \CZ$ is Fréchet differentiable. Moreover, for every point $\xi \in \CXR$, the Fréchet derivative $D \Sol(\xi) \in \CL(L^{\infty}(\Sigma), \CZ)$ is defined as follows: for every direction $\lxi \in L^{\infty}(\Sigma)$, the directional derivative is given by
    \begin{equation*}
        D \Sol(\xi)[\lxi] = (\lv, \lcpm, \lpsi),
    \end{equation*}
    where $(\lv, \lcpm, \lpsi)$ is the solution of the linearized system \eqref{eq:lin_snpp}--\eqref{eq:lin_ic}.
\end{propo}
\begin{proof}
    We want to show that the solution operator to the linearized system \eqref{eq:lin_snpp} is the Fréchet derivative of the control-to-state operator.
    For that, we need to show that the remainder goes to zero faster than linearly, more precisely
    for all $\xi \in \CXR$ and $\lxi \in L^{\infty}(\Sigma)$ we have:
    \begin{align}\label{eq:diff_limit}
        \frac{\norm{\Sol(\xi + \lxi) - \Sol(\xi) - (\lv, \lcpm, \lpsi)}_{\CZ}}{\norm{\lxi}_{L^{\infty}(\Sigma)}} \to 0 \quad \text{ as } \quad \norm{\lxi}_{L^{\infty}(\Sigma)} \to 0,
    \end{align}
    where $(\lv, \lcpm, \lpsi)$ is the solution to the linearized system \eqref{eq:lin_snpp}--\eqref{eq:lin_ic} centered in $(\bv, \cpm, \psi) \coloneqq \Sol(\xi)$ along the direction $\lxi$. Without loss of generality, since $\CXR$ is open and our final aim is to pass to the limit as $\lxi$ goes to $0$, we assume that the norm of $\lxi$ is small enough that $\xi + \lxi \in \CXR$, and
    we denote $\Sol(\xi + \lxi) \eqqcolon  (\bv_{\lxi}, \cpm_{\lxi}, \psi_{\lxi})$. We also introduce the following compact notation
    \begin{equation*}
        \bbv \coloneqq \bvh - \bv -\lv, \qquad \bcpm \coloneqq \cpmh - \cpm - \lcpm, \qquad \bpsi \coloneqq \psih - \psi - \lpsi.
    \end{equation*}
    Then, $(\bbv, \bcpm, \bpsi)$ is a weak solution of the system obtained by subtracting the systems satisfied by $\Sol(\xi)$ and $(\lv, \lcpm, \lpsi)$ from the system satisfied by $\Sol(\xi + \lxi)$. Explicitly,  for the Stokes equation, it holds
        \begin{equation}\label{eq:diff_weak_Stokes}
            \begin{split}
                &\kappa \int_\Omega \partial_t \bbv \cdot \bs{\varphi} \diff \X + \int_\Omega \nabla \bbv : \nabla \bs{\varphi} \diff \X\\
                &\quad = - \int_\Omega \Big[(\ch^+ - \ch^-) \nabla \psih - (c^+ - c^-) \nabla \psi -  (\gamma^+ - \gamma^-) \nabla\psi - ( c^+- c^-)\nabla\lpsi\Big]\cdot \bs{\varphi} \diff \X,\\
                &\quad = - \int_\Omega \Big[(\ch^+ - \ch^-) \nabla \bpsi + (\bc^+ - \bc^-) \nabla \psi  + \Big( (\ch^+ - \ch^-) - ( c^+- c^-) \Big)\nabla\lpsi\Big]\cdot \bs{\varphi} \diff \X,
            \end{split}
        \end{equation}
    for a.e. $t \in (0,T)$, and for every $\bs{\varphi} \in W^{1,2}_{0,\sigma}(\Omega)$, and $\bbv (0) = 0$ a.e. Testing this equation with $\bs{\varphi} = \bbv$, we get:
    \begin{align*}
        &\frac{\kappa}{2} \frac{\diff}{\diff t} \norm{\bbv}_{L^2(\Omega)}^2
        + \norm{\nabla \bbv}_{L^2(\Omega)}^2\\
        &\quad \leq \norm{\ch^+ - \ch^-}_{L^{\infty}(Q)} \norm{\nabla \bpsi}_{L^2(\Omega)} \norm{\bbv}_{L^2(\Omega)} + \sigmp{\norm{\bc^\pm}_{L^2(\Omega)}} \norm{\nabla \psi}_{L^{\infty}(0,T;L^6(\Omega))} \norm{\bbv}_{L^3(\Omega)}\\
        &\quad \quad + \sigmp{\norm{\ch^\pm - \cpm}_{L^{\infty}(0,T;L^2(\Omega))}} \norm{\nabla \lpsi}_{L^{\infty}(0,T;L^6(\Omega))} \norm{\bbv}_{L^3(\Omega)}\\
        &\quad \leq C_R \Big(\norm{\nabla \bpsi}_{L^2(\Omega)} \norm{\bbv}_{L^2(\Omega)} + \sigmp{\norm{\bc^\pm}_{L^2(\Omega)}} \norm{\bbv}_{L^3(\Omega)} + \norm{\lxi}_{L^2(0,T;L^2(\Gamma))} \norm{\lxi}_{L^{\infty}(\Sigma)}  \norm{\bbv}_{L^3(\Omega)}\Big),
    \end{align*}
    where we have bounded the right-hand side by the H\"older inequality and estimates \eqref{eq:main_norm_ineq}--\eqref{eq:cont_dep_ineq}. Then, we employ the Gagliardo--Nirenberg inequality \eqref{eq:GN_L3_0} and the Young inequality with a small constant $\delta > 0$ (fixed but yet to be defined), obtaining:
    \begin{equation} \label{eq:diff_ineq_s}
        \begin{split}
             &\frac{\kappa}{2} \frac{\diff}{\diff t} \norm{\bbv}_{L^2(\Omega)}^2
            + \norm{\nabla \bbv}_{L^2(\Omega)}^2\\
            &\quad \leq \delta \norm{\nabla \bpsi}_{L^2(\Omega)}^2 + C_{\delta,R} \norm{\bbv}_{L^2(\Omega)}^2 +
            C_R \Big(\sigmp{\norm{\bc^\pm}_{L^2(\Omega)}^2} + \norm{\lxi}_{L^2(0,T;L^2(\Gamma))}^2\norm{\lxi}_{L^{\infty}(\Sigma)}^2\\
            &\quad \quad +\norm{\nabla \bbv}_{L^2(\Omega)}\norm{\bbv}_{L^2(\Omega)}\Big)\\
            &\quad \leq \delta \Big(\norm{\nabla \bpsi}_{L^2(\Omega)}^2 + \norm{\nabla \bbv}_{L^2(\Omega)}^2 \Big) + C_{\delta,R} \norm{\bbv}_{L^2(\Omega)}^2\\
            &\quad \quad + C_R \Big(\sigmp{\norm{\bc^\pm}_{L^2(\Omega)}^2}  + \norm{\lxi}_{L^2(0,T;L^2(\Gamma))}^2 \norm{\lxi}_{L^{\infty}(\Sigma)}^2 \Big).
        \end{split}
    \end{equation}
    Regarding the Nernst--Planck equation we get that $\bcpm$ satisfies
    \begin{equation}
        \begin{split}
            \hspace{-1em}&\langle\partial_t\bcpm,\varphi\rangle + \int_\Omega\varMa \nabla \bcpm \cdot \nabla \varphi \diff \X\\
            &\hspace{1em} = \int_\Omega\big( \cpmh \bvh - \cpm \bv -\lcpm \bv - c^{\pm}\lv\big) \cdot \nabla \varphi
            \mp \varMa \big( \cpmh \nabla \psih - \cpm \nabla \psi - \lcpm\nabla\psi -c^{\pm}\nabla\lpsi\big)\cdot\nabla\varphi\diff \X\\
             &\hspace{1em} = \int_\Omega\Big( \bcpm \bv + \cpmh \bbv + (\cpmh - \cpm) \lv \Big) \cdot \nabla \varphi \diff \X \mp \int_\Omega\varMa \Big( \cpmh \nabla \bpsi + \bcpm \nabla \psi + (\cpmh - \cpm ) \nabla\lpsi \Big) \cdot \nabla \varphi \diff \X
        \end{split}
    \end{equation}
    for almost every $t\in(0,T)$ and for all $\varphi\in W^{1,2}(\Omega)$, and $\bcpm(0) = 0$ a.e. We choose $\varphi = \bcpm$ as a test function, and we notice that the first addend on the right-hand side $\int_{\Omega} \bcpm \bv \cdot \nabla \bcpm \diff \X$ vanishes because $\bv$ is divergence-free and $\bv|_\Gamma = 0$. We get:
    \begin{align*}
        &\frac{1}{2}\frac{\diff}{\diff t} \norm{\bcpm}_{L^2(\Omega)}^2
        + \int_\Omega |\nabla \bcpm|^2_{\varMa} \diff \X\\
        &\quad \leq C \Big(\norm{\cpmh}_{L^{\infty}(Q)} \norm{\bbv}_{L^2(\Omega)} + \norm{\cpmh - \cpm}_{L^{\infty}(0,T;L^2(\Omega))} \norm{\lv}_{W^{2,2}(\Omega)} + \norm{\cpmh}_{L^{\infty}(Q)} \norm{\nabla \bpsi}_{L^2(\Omega)}\\
        &\quad\quad+ \norm{\bcpm}_{L^3(\Omega)} \norm{\nabla \psi}_{L^{\infty}(0,T;L^6(\Omega))} + \norm{\cpmh - \cpm}_{L^3(\Omega)} \norm{\nabla \lpsi}_{L^{\infty}(0,T;L^6(\Omega))} \Big) \norm{\nabla \bcpm}_{L^2(\Omega)}\\
        &\quad \leq C \Big(\norm{\bbv}_{L^2(\Omega)} + \norm{\lxi}_{L^2(0,T;L^2(\Gamma))} \norm{\lv}_{W^{2,2}(\Omega)} + \norm{\nabla \bpsi}_{L^2(\Omega)}\\
        &\quad\quad+ \norm{\bcpm}_{L^3(\Omega)} + \norm{\cpmh - \cpm}_{W^{1,2}(\Omega)} \norm{\lxi}_{L^{\infty}(\Sigma)} \Big) \norm{\nabla \bcpm}_{L^2(\Omega)},
    \end{align*}
    having used the H\"older inequality, and the estimates \eqref{eq:main_norm_ineq}--\eqref{eq:cont_dep_ineq}, \eqref{eq:lin_syst_main_est}. 
    Again, the Young inequality with a small $\delta > 0$ and the Gagliardo--Nirenberg inequality \eqref{eq:GN_L3} lead to
    \begin{equation}\label{eq:diff_ineq_np}
        \begin{split}
            &\frac{1}{2}\frac{\diff}{\diff t} \norm{\bcpm}_{L^2(\Omega)}^2
            + \int_\Omega |\nabla \bcpm|^2_{\varMa} \diff \X\\
            &\quad \leq \delta \norm{\nabla \bcpm}_{L^2(\Omega)}^2 + C_{\delta} \Big( \norm{\bbv}_{L^2(\Omega)}^2 + \norm{\lxi}_{L^2(0,T;L^2(\Gamma))}^2 \norm{\lv}_{W^{2,2}(\Omega)}^2 + \norm{\nabla \bpsi}_{L^2(\Omega)}^2\\
            &\quad\quad+ \norm{\bcpm}_{L^2(\Omega)}^2 + \norm{\cpmh - \cpm}_{W^{1,2}(\Omega)}^2 \norm{\lxi}_{L^{\infty}(\Sigma)}^2\Big).
        \end{split}
    \end{equation}
    Regarding the electric potential, $\bpsi$ solves
    \begin{align} \label{eq:diff_weak_p}
        \int_\Omega \varMa \nabla\bpsi \cdot \nabla \varphi \diff \X + \tau \int_\Gamma \bpsi \varphi \dS = \int_\Omega(\bc^+-\bc^-)\varphi\diff \X,
    \end{align}
    for almost every $t\in(0,T)$ and all $\varphi\in W^{1,2}(\Omega)$. 
    Taking $\varphi = \bpsi$, multiplying with a big constant $K > 0$, and bounding the right-hand side by H\"older, generalized Poincaré \eqref{eq:gen_Poincare_ineq}, and Young inequalities, we have
    \begin{equation}\label{eq:diff_ineq_p}
        \begin{split}
            K\int_\Omega |\nabla\bpsi|^2_{\varMa} \diff \X + \tau K\norm{\bpsi}_{L^2(\Gamma)}^2 &\leq K \sigmp{\norm{\bcpm}_{L^2(\Omega)}} \norm{\bpsi}_{L^2(\Omega)}\\
            &\leq K C_P \sigmp{\norm{\bcpm}_{L^2(\Omega)}} \Big( \norm{\nabla \bpsi}_{L^2(\Omega)} + \norm{\bpsi}_{L^2(\Gamma)}\Big)\\
            &\leq \delta K \Big( \norm{\nabla \bpsi}_{L^2(\Omega)}^2 + \norm{\bpsi}_{L^2(\Gamma)}^2\Big) + C_{\delta}K\sigmp{\norm{\bcpm}_{L^2(\Omega)}^2}.
        \end{split}
    \end{equation}
    Then, we sum the inequalities \eqref{eq:diff_ineq_s}, \eqref{eq:diff_ineq_np}, \eqref{eq:diff_ineq_p}, and integrate in time over $(0,t)$ for $t \in (0,T)$, recalling that the initial values vanish. Choosing $\delta$ small enough and, after that, $K$ big enough, it is possible to absorb some of the terms into the corresponding ones on the left-hand side, as we have already done previously. We get:
    \begin{equation*}
        \begin{split}
             & \left( \kappa \norm{\bbv}_{L^2(\Omega)}^2 + \sigmp{\norm{\bcpm}_{L^2(\Omega)}^2} \right) (t)
            + \int_0^t \Big( \norm{\nabla \bbv}_{L^2(\Omega)}^2 + \sigmp{\norm{\nabla \bcpm}_{L^2(\Omega)}^2} + \norm{\nabla \bpsi}_{L^2(\Omega)}^2 + \tau \norm{\bpsi}_{L^2(\Gamma)}^2\Big) \diff s \\
            &\quad \leq C_R\int_0^t \Big( \norm{\bbv}_{L^2(\Omega)}^2 + \sigmp{\norm{\bc^\pm}_{L^2(\Omega)}^2}  \Big) \diff s + C_R t \norm{\lxi}_{L^2(0,T;L^2(\Gamma))}^2 \norm{\lxi}_{L^{\infty}(\Sigma)}^2\\
            &\quad \quad + C_R \norm{\lxi}_{L^2(0,T;L^2(\Gamma))}^2 \int_0^t \norm{\lv}_{W^{2,2}(\Omega)}^2 \diff s + C_R \norm{\lxi}_{L^{\infty}(\Sigma)}^2 \int_0^t \sigmp{\norm{\cpmh - \cpm}_{W^{1,2}(\Omega)}^2} \diff s\\
            &\quad \leq C_R \bigg(\int_0^t \Big( \norm{\bbv}_{L^2(\Omega)}^2 + \sigmp{\norm{\bc^\pm}_{L^2(\Omega)}^2}  \Big) \diff s + \norm{\lxi}_{L^{\infty}(\Sigma)}^2 \norm{\lxi}_{L^2(0,T;L^2(\Gamma))}^2  + \norm{\lxi}_{L^2(0,T;L^2(\Gamma))}^4  \bigg),
        \end{split}
    \end{equation*}
    where the last inequality follows from the estimates \eqref{eq:cont_dep_ineq} and \eqref{eq:lin_syst_main_est}--\eqref{eq:lin_lpsi_add_est}. Applying the Gronwall inequality and the generalized Poincaré inequality \eqref{eq:gen_Poincare_ineq} for $\bpsi$, we get
    \begin{equation}\label{eq:diff_gronwall}
        \begin{split}
            &\norm{\bbv}_{L^{\infty}(0,T;L^2(\Omega)) \cap L^2(0,T;W^{1,2}(\Omega))}^2
            + \sigmp{\norm{\bcpm}_{L^{\infty}(0,T;L^2(\Omega)) \cap L^2(0,T;W^{1,2}(\Omega))}^2}
            + \norm{\bpsi}_{L^2(0,T;W^{1,2}(\Omega))}^2\\
            &\quad \leq C_R \Big( \norm{\lxi}_{L^{\infty}(\Sigma)}^2 \norm{\lxi}_{L^2(0,T;L^2(\Gamma))}^2 + \norm{\lxi}_{L^2(0,T;L^2(\Gamma))}^4 \Big)
            \leq C_R \norm{\lxi}_{L^{\infty}(\Sigma)}^4,
        \end{split}
    \end{equation}
    where in the last inequality we used $\norm{\lxi}_{L^2(0,T;L^2(\Gamma))} \leq |\Sigma|^{1/2} \norm{\lxi}_{L^{\infty}(\Sigma)}$.
    It remains to estimate $\bpsi$ in the stronger norm of $L^{\infty}(0,T;W^{1,6}(\Omega))$ required by $\CZ$.
    By \eqref{eq:diff_weak_p}, for almost every $t \in (0,T)$, the function $\bpsi(t)$ is the weak solution of the Poisson equation $(P_i)$ in \eqref{eq:splitting} with right-hand side given by $\bc^+ - \bc^-$ and homogeneous Robin boundary condition.
    Therefore, the elliptic regularity estimate \eqref{eq:est_psii} with $p = 2$, the continuous embedding $W^{2,2}(\Omega) \hookrightarrow W^{1,6}(\Omega)$, and \eqref{eq:diff_gronwall} yield
    \begin{equation*}
        \norm{\bpsi}_{L^{\infty}(0,T;W^{1,6}(\Omega))}
        \leq C \norm{\bpsi}_{L^{\infty}(0,T;W^{2,2}(\Omega))}
        \leq C \norm{\bc^+ - \bc^-}_{L^{\infty}(0,T;L^2(\Omega))}
        \leq C_R \norm{\lxi}_{L^{\infty}(\Sigma)}^2.
    \end{equation*}
    Combining this estimate with \eqref{eq:diff_gronwall}, we conclude that
    \begin{equation*}
        \norm{(\bbv, \bcpm, \bpsi)}_{\CZ}^2 \leq C_R \norm{\lxi}_{L^{\infty}(\Sigma)}^4,
    \end{equation*}
    from which \eqref{eq:diff_limit} follows.
\end{proof}

\section{First-order optimality conditions}\label{sec:first_order}

In this section, we conclude the derivation of the necessary first-order conditions characterizing optimal controls. Going back to the inequality \eqref{eq:opt_cond_0}---which is now rigorous---we compute the derivative using the chain rule, obtaining the following result.

\begin{propo} \label{propo:opt_cond_0}
    Let $\xi^* \in \CXad$ be an optimal control for \eqref{eq:opt_control_prob} with corresponding state $ (\bv^*, {\cpm}^*, \psi^*)$.
    Then, the necessary first-order optimality conditions hold
    \begin{equation}\label{eq:opt_cond_1}
        \begin{split}
             D\redJ(\xi^*)[\xi - \xi^*] = &- \alpv  \int_0^T \int_\Omega (\nabla \times \bv^*) \cdot (\nabla \times \lv) \diff \X \diff t \\
            &\quad+ \int_0^T \int_\Omega \sigmp{ \alpcpm \big({\cpm}^* - \cpm_Q\big) \lcpm} \diff \X \diff t   + \alpxi \int_0^T \int_\Gamma \xi^* (\xi - \xi^*) \dS \diff t \geq 0 
        \end{split}
    \end{equation}
    for all $\xi \in \CXad$, where $(\lv, \lcpm, \lpsi) \coloneqq D\Sol(\xi^*)[\xi - \xi^*]$.
\end{propo}

\begin{remark}
     These are only necessary conditions, as we do not have convexity of the reduced cost functional: even if $\alpv  = 0$ and  $\J$ is convex, the solution operator $\Sol$ is highly nonlinear and, as a result, the composition $\redJ$ is in general nonconvex, compromising the application of stronger results.
\end{remark}

\noindent Unfortunately, this characterization is not completely satisfactory, since practically computing the derivative $D\Sol(\xi^*)[\xi - \xi^*]$ appearing in \eqref{eq:opt_cond_1} requires solving an infinite number of linearized systems \eqref{eq:lin_snpp}--\eqref{eq:lin_ic}, one for each $\xi \in \CXad$. 
Thus, our next goal is to rewrite \eqref{eq:opt_cond_1} in an equivalent way, 
removing the linearized variables. This is classically done in the literature by means of the \textit{adjoint system}, which can be formally derived using the Lagrangian method (see, e.g., \cite{Troltzsch_2010}). 
More explicitly, writing the Lagrangian functional as the sum of the cost functional and the SNPP system \eqref{eq:snpp}--\eqref{eq:ic} multiplied by the adjoint variables $(\av, \ap, \acpm, \apsi)$ (which play the role of Lagrangian multipliers), and setting its derivative with respect to the state variables to zero, one obtains the following backward-in-time system: 
\begin{subequations}\label{eq:adjoint_snpp}
    \begin{align}
        - \kappa \,  \partial_t \av - \Delta \av + \nabla \ap - \sigmp{\cpm \, \nabla \acpm} &= - \alpv \nabla \times (\nabla \times \bv) &&\text{ in } Q,\label{eq:adjoint_s}\\
        \nabla \cdot \av &= 0 &&\text{ in } Q, \label{eq:adjoint_incompress}\\
        - \partial_t \acpm - \bv \cdot \nabla \acpm - \nabla \cdot (\varMa \nabla \acpm)
        \pm \varMa \nabla \psi \cdot \nabla \acpm \pm \nabla \psi \cdot \av \mp \apsi &= \alpcpm ({\cpm} - \cpm_Q)&&\text{ in } Q, \label{eq:adjoint_np}\\
        - \nabla \cdot (\varMa \nabla \apsi) + \sigmp{\mp \nabla \cdot (\cpm \varMa \nabla \acpm)} - \nabla \cdot ((c^+ - c^-) \av) &= 0 &&\text{ in } Q, \label{eq:adjoint_p}
\intertext{equipped with the boundary conditions\newsubeqgroup{eq:adjoint_bc}}
        \av &= 0 &&\text{ on } \Sigma,\\ 
        \varMa \nabla \acpm \cdot \bn &= 0 &&\text{ on } \Sigma,\\
        \varMa \nabla \apsi \cdot \bn + \tau \apsi &= 0  &&\text{ on } \Sigma, 
\intertext{and the final conditions \stopsubeqgroup}
        \av(T) = 0, \quad \acpm(T) &= 0 &&\text{ in } \Omega. \label{eq:adjoint_fc}
    \end{align}
\end{subequations}

The system \eqref{eq:adjoint_snpp}--\eqref{eq:adjoint_fc} has to be interpreted in a weak sense, according to the following definition.

\begin{definition}[Weak solution of the adjoint system]
\label{definition:adjoint_weak_sol}
    Let $(\bv, \cpm, \psi) \in \CY$.  A tuple $(\av,\acpm,\apsi)$ with the regularity
    \begin{align*}
        \av &\in W^{1,2}(0,T;L^2(\Omega)) \cap L^2(0,T;W^{1,2}_{0,\sigma}(\Omega) \cap W^{2,2}(\Omega)),\\
        \acpm &\in W^{1,2}(0,T;(W^{1,2}(\Omega))^*) \cap C^0([0,T];L^2(\Omega)) \cap L^2(0,T;W^{1,2}(\Omega)),\\
        \apsi &\in L^2(0,T;W^{1,2}(\Omega))
    \end{align*}
    is called a weak solution of \eqref{eq:adjoint_snpp}--\eqref{eq:adjoint_fc} if the system \eqref{eq:adjoint_snpp}--\eqref{eq:adjoint_bc} is satisfied in the following variational sense
    \begin{subequations} \label{eq:adjoint_weak}
        \begin{align}
          &- \kappa\,\langle\partial_t\av,\bs{\varphi}\rangle + \int_\Omega\nabla\av:\nabla\bs{\varphi}\diff \X - \int_\Omega \sigmp{\cpm \, \nabla \acpm} \cdot\bs{\varphi}\diff \X = - \alpv  \int_\Omega (\nabla \times \bv) \cdot (\nabla \times \bs{\varphi}) \diff \X, \label{eq:adjoint_weak_s}\\
         &\begin{aligned}
              &-\langle\partial_t\acpm,\varphi\rangle - \int_\Omega \bv \cdot \nabla \acpm \varphi \diff \X + \int_\Omega\varMa \nabla \acpm \cdot \nabla \varphi \diff \X
              \pm \int_\Omega\big( \varMa \nabla \psi \cdot \nabla \acpm + \nabla \psi \cdot \av - \apsi \big) \varphi \diff \X\\
              &\quad \quad =  \alpcpm \int_\Omega({\cpm} - \cpm_Q) \varphi \diff \X,
         \end{aligned} \label{eq:adjoint_weak_np}\\
          &\int_\Omega \varMa \nabla\apsi\cdot\nabla\varphi\diff \X + \tau \int_\Gamma \apsi \, \varphi \dS + \int_\Omega \sigmp{\pm \cpm \varMa \nabla \acpm} \cdot \nabla \varphi \diff \X + \int_\Omega(c^+-c^-) \av \cdot \nabla \varphi\diff \X = 0, \label{eq:adjoint_weak_p}
        \end{align}
    \end{subequations}
    for almost every $t\in(0,T)$ and for all $\bs{\varphi}\in W^{1,2}_{0,\sigma}(\Omega)$, $\varphi \in W^{1,2}(\Omega)$, and if the final conditions \eqref{eq:adjoint_fc} hold in $L^2(\Omega)$.
\end{definition}

\begin{remark}
    Once again, the adjoint pressure $\ap$ can be a posteriori reconstructed up to an additive constant, as it was for the state and the linearized cases (cf.~\zcref{remark:p} and \zcref{remark:pi}). 
    Choosing $\ap$ with $0$ mean value, by comparison in \eqref{eq:adjoint_s}, it follows that $\ap \in L^2(0,T;W^{1,2}(\Omega))$.
\end{remark}

\begin{propo}
    For every $\xi \in \CXR$ with corresponding state $\Sol(\xi) \eqqcolon (\bv, \cpm, \psi)$, the adjoint system \eqref{eq:adjoint_snpp}--\eqref{eq:adjoint_fc} has a unique weak solution in the sense of \zcref{definition:adjoint_weak_sol}.
\end{propo}
\begin{proof}
    The existence of a weak solution then follows via a standard Galerkin approximation and weak compactness arguments. Since it is quite standard, here we only derive formal a priori estimates. 

    \textbf{Poisson estimate.} 
    We test \eqref{eq:adjoint_weak_p} with $\apsi$ and estimate the right-hand side using the boundedness of $\varMa$, the estimate \eqref{eq:bound_charges}, the H\"older and the Young inequalities. We obtain:
    \begin{equation*}
        \begin{split}
             \int_\Omega |\nabla \apsi|^2_{\varMa} \diff \X + \tau \norm{\apsi}_{L^2(\Gamma)}^2  &=  - \int_\Omega \sigmp{ \pm \cpm \varMa \nabla \acpm} \cdot \nabla \apsi \diff \X
            - \int_\Omega (c^+ - c^-) \av \cdot \nabla \apsi \diff \X\\
            & \leq C \Big(\sigmp{\norm{\cpm}_{L^{\infty}(Q)} \norm{\nabla \acpm}_{L^2(\Omega)}} + \sigmp{\norm{\cpm}_{L^{\infty}(Q)}} \norm{\av}_{L^2(\Omega)} \Big) \norm{\nabla \apsi}_{L^2(\Omega)}\\
            & \leq C_R \Big(\sigmp{\norm{\nabla \acpm}_{L^2(\Omega)}} + \norm{\av}_{L^2(\Omega)} \Big) \norm{\nabla \apsi}_{L^2(\Omega)}\\
            & \leq \delta \norm{\nabla \apsi}_{L^2(\Omega)}^2 + C_{R, \delta} \Big(\sigmp{\norm{\nabla \acpm}^2} + \norm{\av}_{L^2(\Omega)}^2 \Big)
        \end{split} 
    \end{equation*}
    for a small constant $\delta > 0$ to be fixed small enough so that the addend multiplied by it can be absorbed into the left-hand side. Exploiting the Poincar\'{e} inequality \eqref{eq:gen_Poincare_ineq}, we get
    \begin{align}\label{eq:adj_poisson_est}
        \norm{\apsi}_{W^{1,2}(\Omega)}^2 \leq C_R \Big(\sigmp{\norm{\nabla \acpm}_{L^2(\Omega)}^2} + \norm{\av}_{L^2(\Omega)}^2 \Big)
    \end{align}
    for a.e. $t \in (0,T)$.
    
    \paragraph{Stokes estimate.}
    Testing \eqref{eq:adjoint_weak_s} with $\av$ and using H\"{o}lder's inequality, Young's inequality with a small $\delta > 0$, the boundedness of the charges \eqref{eq:bound_charges} and estimate \eqref{eq:main_norm_ineq} give
    \begin{equation}
        \begin{split}\label{eq:adj_stokes_est}
            - \frac{1}{2} &\frac{\diff}{\diff t} \norm{\av}_{L^2(\Omega)}^2
            + \norm{\nabla \av}_{L^2(\Omega)}^2 
            = - \alpv \int_\Omega (\nabla \times \bv) \cdot (\nabla \times \av) \diff \X 
                + \int_\Omega \sigmp{\cpm \, \nabla \acpm \cdot \av} \diff \X\\
            &\quad \leq C \norm{\nabla \bv}_{L^2(\Omega)} \norm{\nabla \av}_{L^2(\Omega)} + \sigmp{\norm{\cpm}_{L^{\infty}(Q)} \norm{\nabla \acpm}_{L^2(\Omega)}} \norm{\av}_{L^2(\Omega)}\\
            &\quad \leq C_R \Big( \norm{\nabla \bv}_{L^2(\Omega)}\norm{\nabla \av}_{L^2(\Omega)} + \sigmp{\norm{\nabla \acpm}_{L^2(\Omega)}} \norm{\av}_{L^2(\Omega)} \Big)\\
            &\quad \leq \delta \Big( \norm{\nabla \av}_{L^2(\Omega)}^2 + \sigmp{\norm{\nabla \acpm}_{L^2(\Omega)}^2} \Big) + C_{R,\delta} \Big( \norm{\nabla \bv}_{L^2(\Omega)}^2 + \norm{\av}_{L^2(\Omega)}^2 \Big)
        \end{split}
    \end{equation}

   \paragraph{Nernst--Planck estimate.} 
    We test \eqref{eq:adjoint_weak_np} with $\acpm$ for each species and sum over $\pm$. The convection terms vanish since $\bv$ is divergence-free and $\bv = 0$ on $\Gamma$. Estimating the right-hand side using the boundedness of $\varMa$ by assumption \ref{ass:anisotropy} and the H\"older inequality, we get:
    \begin{equation*}
        \begin{split}
            & - \frac{1}{2} \frac{\diff}{\diff t} \sigmp{\norm{\acpm}_{L^2(\Omega)}^2}
            + \int_\Omega\sigmp{|\nabla \acpm|_{\varMa}^2} \diff \X\\
            &\quad = \int_\Omega 
            \sigmp{ \mp \varMa \nabla \psi \cdot \nabla \acpm \, \acpm \mp (\nabla \psi \cdot \av) \acpm \pm \apsi \,\acpm} \diff \X+ \int_\Omega \sigmp{\alpcpm ({\cpm} - \cpm_Q) \acpm} \diff \X\\
            &\quad \leq C \bigg( \norm{\nabla \psi}_{L^{\infty}(0,T;L^6(\Omega))} \sigmp{\norm{\nabla \acpm}_{L^2(\Omega)} \norm{\acpm}_{L^3(\Omega)}} + \norm{\nabla \psi}_{L^{\infty}(0,T;L^6(\Omega))} \norm{\av}_{L^2(\Omega)} \sigmp{\norm{\acpm}_{L^3(\Omega)}}\\
            &\quad\quad + \norm{\apsi}_{L^2(\Omega)} \sigmp{\norm{\acpm}_{L^2(\Omega)}} + \sigmp{\norm{\cpm -\cpm_Q}_{L^2(\Omega)}\norm{\acpm}_{L^2(\Omega)}}\bigg)
        \end{split}
    \end{equation*}
    Then, we exploit the boundedness of $\cpm_Q$ in $L^2(Q)$ by assumption \ref{ass:cost_funct}, estimates \eqref{eq:main_norm_ineq}--\eqref{eq:bound_charges}, the Gagliardo--Nirenberg inequality \eqref{eq:GN_L3}, and the Young inequality with a small constant $\delta>0$. We obtain:
    \begin{equation}\label{eq:adj_charges_est}
        \begin{split}
            & - \frac{1}{2} \frac{\diff}{\diff t} \sigmp{\norm{\acpm}_{L^2(\Omega)}^2}
            + \int_\Omega\sigmp{|\nabla \acpm|_{\varMa}^2} \diff \X\\
            &\quad \leq C_R \Big(\sigmp{\norm{\nabla \acpm}_{L^2(\Omega)} \norm{\acpm}_{L^3(\Omega)}} + \norm{\av}_{L^2(\Omega)} \sigmp{\norm{\acpm}_{L^3(\Omega)}}\\
            &\quad\quad + \norm{\apsi}_{L^2(\Omega)} \sigmp{\norm{\acpm}_{L^2(\Omega)}} + \sigmp{\norm{\cpm -\cpm_Q}_{L^2(\Omega)}\norm{\acpm}_{L^2(\Omega)}}\Big)\\
            &\quad \leq \delta \Big(\sigmp{\norm{\nabla \acpm}_{L^2(\Omega)}^2} +\norm{\apsi}_{L^2(\Omega)}^2 \Big)  + C_{R,\delta} \Big(\sigmp{\norm{\acpm}_{L^3(\Omega)}^2} + \sigmp{\norm{\acpm}_{L^2(\Omega)}^2} \Big)\\
            &\quad \quad + C \Big(\norm{\av}_{L^2(\Omega)}^2  + \sigmp{\norm{\cpm -\cpm_Q}_{L^2(\Omega)}^2} \Big)\\
            &\quad \leq \delta \Big(\sigmp{\norm{\nabla \acpm}_{L^2(\Omega)}^2} +\norm{\apsi}_{L^2(\Omega)}^2 \Big) + C_{R,\delta} \sigmp{\norm{\acpm}_{L^2(\Omega)}^2} + C \Big(\norm{\av}_{L^2(\Omega)}^2 + \sigmp{\norm{\cpm -\cpm_Q}_{L^2(\Omega)}^2} \Big).
        \end{split}
    \end{equation}
    We sum the inequalities \eqref{eq:adj_poisson_est} multiplied by a small constant $\tilde{\delta}>0$ and \eqref{eq:adj_stokes_est}--\eqref{eq:adj_charges_est} integrated in time over the interval $(t,T)$ for a $t \in (0,T)$. 
    Then, we choose $\tilde{\delta}$ small enough (so that $\tilde{\delta}C_R\sigmp{\norm{\nabla \acpm}_{L^2(t,T;L^2(\Omega))}^2}$ can be absorbed into the corresponding term on the left-hand side) and, subsequently, $\delta$ small enough so that the terms multiplied by it on the right-hand side can be absorbed into the corresponding terms on the left-hand side.
    
   \paragraph{Gronwall estimate.}
    This way we obtain: 
    \begin{equation*}
        \begin{split}
            &\left(\norm{\av}_{L^2(\Omega)}^2 + \sigmp{\norm{\acpm}_{L^2(\Omega)}^2}\right)(t) + \int_t^T \Big(\norm{\apsi}_{W^{1,2}(\Omega)}^2 + \norm{\nabla \av}_{L^2(\Omega)}^2 \Big) \diff t + \int_t^T \int_\Omega\sigmp{|\nabla \acpm|_{\varMa}^2} \diff \X \diff t\\
            &\quad \leq C_R \int_t^T \Big(\norm{\av}_{L^2(\Omega)}^2 +\sigmp{\norm{\acpm}_{L^2(\Omega)}^2} \Big) \diff t + C_R \int_t^T \Big( \norm{\nabla \bv}_{L^2(\Omega)}^2 + \sigmp{\norm{\cpm -\cpm_Q}_{L^2(\Omega)}^2}\Big) \diff t,
        \end{split}
    \end{equation*}
    where we also exploited the vanishing final conditions \eqref{eq:adjoint_fc}.
    Applying Gronwall's inequality yields
    \begin{equation}\label{eq:adj_main_est}
       \begin{split}
            &\norm{\apsi}_{L^2(0,T;W^{1,2}(\Omega))}^2 + \norm{\av}_{L^{\infty}(0,T;L^2(\Omega)) \cap L^2(0,T;W^{1,2}(\Omega))}^2 + \sigmp{\norm{\acpm}_{L^{\infty}(0,T;L^2(\Omega)) \cap L^2(0,T;W^{1,2}(\Omega))}^2}\\
            &\quad \leq C_R \Big(\norm{\nabla \bv}_{L^2(0,T;L^2(\Omega))}^2 + \sigmp{\norm{\cpm -\cpm_Q}_{L^2(0,T;L^2(\Omega))}^2}\Big) e^{C_RT}.
       \end{split}
    \end{equation}

    \paragraph{Additional regularity.} 
    Regarding the adjoint Stokes equation \eqref{eq:adjoint_s}--\eqref{eq:adjoint_incompress}, we notice that $\sigmp{\cpm \nabla \acpm}$ is uniformly bounded in $L^2(0,T;L^2(\Omega))$ thanks to \eqref{eq:bound_charges} and \eqref{eq:adj_main_est}, and that $\nabla \times (\nabla \times \bv)$ is uniformly bounded in $L^2(0,T;L^2(\Omega))$ thanks to \eqref{eq:main_norm_ineq}. Thus, applying \cite[Section IV, Thm.~2.7.1]{sohr_01} reverse in time (i.e., for $\tilde{\av}(t, \cdot) \coloneqq \av(T-t,\cdot)$ and the corresponding equation), we have
    \begin{equation}
        \begin{split}
            &\norm{\partial_t \av}_{L^2(0,T;L^2(\Omega))} + \norm{\av}_{ L^2(0,T;W^{2,2}(\Omega))}\\
            &\quad \leq C \bigg( \sigmp{\norm{\cpm \nabla \acpm}_{L^2(0,T;L^2(\Omega))}} + \norm{\nabla \times (\nabla \times \bv)}_{L^2(0,T;L^2(\Omega))} \bigg) \leq C_R
        \end{split}
    \end{equation}
    because $\av(T) = 0$ by \eqref{eq:adjoint_fc}. Turning our attention to the adjoint Nernst--Planck equation \eqref{eq:adjoint_np}, we rewrite it as 
    \begin{equation}\label{eq:adj_np_2}
        - \partial_t \acpm  = \nabla \cdot (\varMa \nabla \acpm) + \alpcpm (\cpm - \cpm_Q) + (\bv \mp \varMa \nabla \psi) \cdot \nabla \acpm  \mp \nabla \psi \cdot \av  \pm \apsi.
    \end{equation}
    First, we notice that  $\varMa \nabla \acpm$ is uniformly bounded in $L^2(0,T;L^2(\Omega))$ by \eqref{eq:adj_main_est} and hypotheses \ref{ass:anisotropy}. 
    Thus, the first addend on the right-hand side is uniformly bounded in $L^2(0,T;(W^{1,2}(\Omega))^*)$. 
    The second term and the last term are uniformly bounded in $L^2(0,T;L^2(\Omega)) \hookrightarrow L^2(0,T;(W^{1,2}(\Omega))^*)$ by \eqref{eq:main_norm_ineq}, assumption \ref{ass:cost_funct}, and \eqref{eq:adj_main_est}.
    Moreover, we observe that, since $W^{1,2}(\Omega) \hookrightarrow L^p(\Omega)$ for $p \in [1,6]$, then $L^{p'}(\Omega) \hookrightarrow (W^{1,2}(\Omega))^*$ for $p' \in [6/5,+\infty]$, in particular for $p' = 3/2$. Thus, we have
    \begin{equation*}
        \begin{split}
            &\norm{(\bv \mp \varMa \nabla \psi) \cdot \nabla \acpm  \mp \nabla \psi \cdot \av}_{L^2(0,T;(W^{1,2}(\Omega))^*)} \leq C \norm{(\bv \mp \varMa \nabla \psi) \cdot \nabla \acpm  \mp \nabla \psi \cdot \av}_{L^2(0,T;L^{3/2}(\Omega))}\\
            &\quad \leq C \bigg(\Big(\norm{\bv}_{L^{\infty}(0,T;L^6(\Omega))} + \norm{\nabla \psi}_{L^{\infty}(0,T;L^6(\Omega))} \Big)\norm{\nabla \acpm}_{L^2(0,T;L^2(\Omega))}\\
            &\quad\quad + \norm{\nabla \psi}_{L^{\infty}(0,T;L^6(\Omega))} \norm{\av}_{L^2(0,T;L^2(\Omega))}\bigg) \leq C_R,
        \end{split}
    \end{equation*}
    having used the H\"older inequality for $1/6 + 1/2 = 2/3$, the estimates \eqref{eq:main_norm_ineq} and \eqref{eq:adj_main_est}, and the continuous embedding
    \begin{equation*}
        W^{1,2}(0,T;L^2(\Omega)) \cap L^2(0,T;W^{2,2}(\Omega)) \hookrightarrow C^0([0,T];W^{1,2}(\Omega)) \hookrightarrow C^0([0,T];L^6(\Omega))
    \end{equation*}
    for the velocity $\bv$.
    By comparison in \eqref{eq:adj_np_2}, it follows
    \begin{equation}
        \norm{\partial_t \acpm}_{L^2(0,T;(W^{1,2}(\Omega))^*)} \leq C_R.
    \end{equation}
    
    \paragraph{Uniqueness.} Uniqueness follows from the linearity of the system and the estimate \eqref{eq:adj_main_est}. More explicitly, considering two solutions, their difference satisfies \eqref{eq:adjoint_snpp}--\eqref{eq:adjoint_fc} with the right-hand sides of equations \eqref{eq:adjoint_s} and \eqref{eq:adjoint_np} equal to zero. Thus, the estimate \eqref{eq:adj_main_est} holds with a vanishing right-hand side, and the two solutions must coincide. 
\end{proof}

\noindent We can now state the first-order optimality conditions in terms of the adjoint variables.

\begin{thm}[First-order optimality conditions---adjoint formulation]\label{thm:opt_cond}
    Let $\xi^* \in \CXad$ be an optimal control for the minimization problem \eqref{eq:opt_control_prob} with corresponding state $ (\bv^*, {\cpm}^*, \psi^*) \coloneqq \Sol(\xi^*)$.
    Then, for every $\xi \in \CXad$, it holds 
    \begin{align} \label{eq:opt_cond}
        D\redJ(\xi^*)[\xi - \xi^*] = \int_0^T \int_\Gamma (\alpxi \xi^* + \apsi) (\xi - \xi^*) \dS \diff t \geq 0,
    \end{align}
    where $\apsi$ is the last component of the solution of the adjoint problem \eqref{eq:adjoint_snpp}--\eqref{eq:adjoint_fc} written for $\xi = \xi^*$, $(\bv, \cpm, \psi) = (\bv^*, {\cpm}^*, \psi^*)$.
\end{thm}

\begin{remark}
    This is exactly the optimality condition we would obtain by formally differentiating the Lagrangian introduced above with respect to the control $\xi$.
\end{remark}

\begin{remark}
    If $\alpxi > 0$, another way to interpret \zcref{thm:opt_cond} is saying that every optimal control is the $L^2(0,T;L^2(\Gamma))$-projection of the corresponding adjoint variable $-\apsi/\alpxi$ on the closed (in $L^2(0,T;L^2(\Gamma))$) and convex set $\CXad$ (see, e.g., \cite[Theorem 5.2]{brezis_11}). Moreover, this projection can be computed pointwise, i.e.,
    \begin{equation*}
        \xi^*(t,x) = \mathbb{P}_{[\lbxi(t,x),\ubxi(t,x)]} \Big( - \frac{1}{\alpxi}\apsi(t,x)\Big) \qquad \text{for a.e. } (t,x) \in \Sigma
    \end{equation*}
    where, for any real numbers $r_1 \leq r_2$, the projection of $\R$ onto the interval $[r_1,r_2]$ is defined as $\mathbb{P}_{[r_1,r_2]}(r) \coloneqq \min \big\{ r_2, \max\{r_1, r\}\big\}$ (for more details, see \cite[Theorem 2.28]{Troltzsch_2010}). More generally, if $\alpxi \geq 0$, it holds
    \begin{equation*}
        \xi^*(t,x) = \begin{cases}
            \lbxi(t,x) &\quad \text{ if } \alpxi \xi^*(t,x) + \apsi(t,x) > 0,\\
            \in [\lbxi(t,x), \ubxi(t,x)] &\quad \text{ if } \alpxi \xi^*(t,x) + \apsi(t,x) = 0,\\
            \ubxi(t,x) &\quad \text{ if } \alpxi \xi^*(t,x) + \apsi(t,x) < 0\\
        \end{cases}
    \end{equation*}
    (see \cite[Lemma 2.26]{Troltzsch_2010}). In particular, if $\alpxi = 0$ and $\apsi \neq 0$ a.e., then $\xi^*$ is a so-called \textit{bang-bang control}.
\end{remark}

\begin{proof}
    To prove that \eqref{eq:opt_cond} holds, it is enough to start from \zcref{propo:opt_cond_0} and show that
    \begin{equation}\label{eq:derivative_id}
        - \alpv \int_{0}^{T} \int_\Omega (\nabla \times \bv^*) \cdot (\nabla \times \lv) \diff \X \diff t + \int_0^T \int_\Omega  \sigmp{\alpcpm({\cpm}^* - \cpm_Q)\lcpm} \diff \X \diff t 
        = \int_0^T \int_\Gamma  (\xi - \xi^*) \apsi \dS \diff t.
    \end{equation}
    To derive this identity, we test the linearized system \eqref{eq:lin_weak} (centered in $\xi^*$ and with $h = \xi - \xi^*$) with the adjoint variables $(\av, \acpm, \apsi)$ and we subtract the adjoint system \eqref{eq:adjoint_weak} (written for $\xi^*$) tested with the solution to the linearized system $(\lv, \lcpm, \lpsi)$. After canceling terms, we obtain: 
    \begin{equation*}
        \begin{split}
            &\frac{\diff}{\diff t} \int_\Omega \Big( \kappa \lv \av + \sigmp{\lcpm \acpm} \Big) \diff \X\\
            &\quad = \int_\Gamma (\xi - \xi^*) \apsi \dS + \int_\Omega \alpv (\nabla \times \bv^*) \cdot (\nabla \times \lv) \diff \X -   \int_\Omega \sigmp{\alpcpm ({\cpm}^* - \cpm_Q) \lcpm} \diff \X.
        \end{split}
    \end{equation*}
    Integrating in time over $(0,T)$, the left-hand side vanishes thanks to the initial conditions \eqref{eq:lin_ic} and to the final conditions \eqref{eq:adjoint_fc}. Thus, we have proved \eqref{eq:derivative_id} and the proof is finished. 
\end{proof}

\section*{Declaration of generative AI use}
During the preparation of this work, the authors used various Claude (Anthropic) models (up to Opus 5.5) for proofreading and language editing, for double-checking selected computations, and for suggesting the proof strategy for the boundedness of the charges in Step~\ref{enum:bound_charges} of the proof of \zcref{thm:weak_existence}.
The authors reviewed all AI-generated suggestions and carried out all proof details.
The authors take full responsibility for the content and correctness of this publication.

\section*{Acknowledgments}
The authors would like to thank Robert Lasarzik for several useful discussions and suggestions which have significantly improved this paper. %
Both the authors are funded by the Deutsche Forschungsgemeinschaft (DFG, German Research Foundation) under Germany's Excellence Strategy – The Berlin Mathematics Research Center MATH+ (EXC-2046/2, project ID: 390685689).
 
\printbibliography
    
\end{document}